\documentclass[11pt,reqno]{amsart}
\usepackage{amsmath,amssymb,amsthm,mathtools}
\usepackage[a4paper,margin=29mm]{geometry}
\usepackage{microtype,array,float}
\usepackage[numbers,sort&compress]{natbib}
\usepackage{hyperref}
\usepackage{xurl}
\hypersetup{colorlinks=true,linkcolor=blue,citecolor=blue,urlcolor=blue,
 pdftitle={A Goldbach-Type Problem for Random Multiplicative Functions}}
\numberwithin{equation}{section}
\newtheorem{theorem}{Theorem}[section]
\newtheorem{lemma}[theorem]{Lemma}
\newtheorem{corollary}[theorem]{Corollary}
\newtheorem{proposition}[theorem]{Proposition}
\theoremstyle{definition}
\newtheorem{definition}[theorem]{Definition}
\newtheorem{question}{Question}[section]
\theoremstyle{remark}
\newtheorem{remark}{Remark}[section]
\newtheorem*{identitylemma}{Basic identities}
\newcommand{\PP}{\mathbb P}
\newcommand{\EE}{\mathbb E}
\newcommand{\ind}{\mathbf 1}
\newcommand{\eps}{\varepsilon}

\title[Goldbach-type problem for random multiplicative functions]{A Goldbach-Type Problem for Random Multiplicative Functions}
\author{Tian Cao Lin}
\address{Institute of Computational Mathematics and Scientific/Engineering Computing,
Academy of Mathematics and Systems Science, Chinese Academy of Sciences;School of Mathematical Sciences, University of Chinese Academy of Sciences}
\email{tiancaolin@lsec.cc.ac.cn}
\date{November 22, 2026}
\begin{document}
\begin{abstract}
We study the counts of the four sign patterns of a random multiplicative function on ordered decompositions $a+b=N$ into positive integers.
The main term is $N/4$ in the nonvanishing model and $F(N)/4$ in the squarefree-supported model,
where $F(N)$ counts the pairs with both summands squarefree and has asymptotic main term $\mathfrak S_{\mathrm{sf}}(N)N$.
Let $\rho_q$ be the probability of a negative sign at the prime $q$. Set $m_q=\min(\rho_q,1-\rho_q)$,
$A_N=\sum_{(\log N)^3<q\le N}m_q/q$, and $B(N)=\sum_{q\le N}m_q/q$.
In the general model with independent prime vectors, almost surely for all sufficiently large integers $N$,
the error for opposite-sign pairs is $O_{\mathcal X}(Ne^{-4A_N})$,
and the error for same-sign pairs is $O_{\mathcal X}(N(e^{-4A_N}+e^{-2B(N)}))$.
In four specified models, if the biases of the independent basic signs are square summable, then for every $\varepsilon>0$
the error is almost surely $O_{\varepsilon,\mathcal X}(N^{1/2+\varepsilon})$.
In the completely multiplicative model under the same bias condition, we prove that
$C_f(N)/\sqrt{N\log N}$ almost surely has limit superior $2$ and limit inferior $-2$.
Together with a known bound for the one-point sum, this gives limit superior $1/2$
and limit inferior $-1/2$ for each normalized counting error about $N/4$.
\end{abstract}
\maketitle
\setcounter{tocdepth}{1}
\tableofcontents

\section{Introduction}
\subsection{The representation problem}
Given a random multiplicative function $f$, we study
\begin{equation}\label{eq:R}
R^f_{\sigma\tau}(N)=\sum_{a=1}^{N-1}
\ind_{\{f(a)=\sigma\}}\ind_{\{f(N-a)=\tau\}},
\qquad \sigma,\tau\in\{-1,1\}.
\end{equation}
We count ordered decompositions $N=a+b$ in which $a,b$ are positive integers and may be equal.
The classical Goldbach problem asks whether both summands are prime. Here we consider their function values, for example
$f(a)=+1$ and $f(b)=-1$.

\noindent\textbf{Motivation.}
The probability of the value $-1$ at a prime plays a central role. Denote this probability by $\rho$ and put $\theta=1-2\rho$.
A squarefree integer with $k$ prime factors then has negative sign with probability $(1-\theta^k)/2$, which in general differs from $\rho$.
For squarefree $a,b$, put $d=(a,b)$. The probability that both signs are negative is
\[
\frac{1-\theta^{\omega(a)}-\theta^{\omega(b)}
+\theta^{\omega(a)+\omega(b)-2\omega(d)}}4.
\]
We take every zeroth power to be $1$. For squarefree $n$, write $g(n)=\prod_{q\mid n}\xi_q$.
In $g(a)g(b)$ every common prime factor contributes its sign twice, giving a product of $1$. Thus
$\EE g(a)=\theta^{\omega(a)}$, $\EE g(b)=\theta^{\omega(b)}$,
and $\EE[g(a)g(b)]=\theta^{\omega(a)+\omega(b)-2\omega(d)}$.
Expanding the two indicators for a negative sign gives the displayed probability formula.
When $0<\rho<1$, we have $|\theta|<1$. If each summand has a growing number of prime factors that do not divide the other,
all three powers tend to zero and the pair probability tends to $1/4$. We therefore need estimates for the mean sign at a single integer and for the correlation between the two summands.

\subsection{Related work}
\noindent\textbf{Random multiplicative sums and sign patterns.}
The study of random multiplicative functions goes back to Wintner \cite{Wintner}.
K\'atai \cite{Katai} studied the normality of random completely multiplicative functions and pointed out that one-point balance is insufficient for uniform pattern frequencies.
Najnudel \cite{Najnudel} studied joint distributions of consecutive values and almost-sure empirical distributions.
Aymone and Sidoravicius \cite{AymoneSidoravicius} considered partial sums in a biased squarefree model.
Harper \cite{Harper} used martingale and moment methods to analyze random multiplicative sums grouped by the number of prime factors.
Klurman, Shkredov and Xu \cite{KSX} and Hardy \cite{Hardy} studied, respectively,
Steinhaus random multiplicative sums at fixed polynomials and exponential sums with random multiplicative coefficients.
In the reflection family $(f(a),f(N-a))$ considered here, the second linear form varies with the target $N$.
This differs from the fixed-length consecutive patterns $(f(n+1),\ldots,f(n+k))$ studied in related settings.

\medskip\noindent\textbf{Goldbach-type additive problems.}
Pylaiev \cite[Theorem 1.13 and Corollary 1.14]{Pylaiev} proved a Gaussian limit for the reflection sum in the fair squarefree model
and a distributional result at square-root scale for opposite-sign representation counts.
De la Bret\`eche, Wang, and Xu \cite[Theorems~1.2 and~7.1]{dBWXFixedPolynomial}
proved Gaussian limits for fixed polynomial sums in the fair completely multiplicative
model, including those for separable quadratic polynomials with positive leading coefficient.
Here the reflection polynomial $a(N-a)$ varies with the target $N$.
Theorem~\ref{thm:cm-extrema} controls distinct targets on the same infinite realization
and determines both signed extremes and their short-interval limits.
Tian and Yang \cite{TianYang} studied sign representations for the deterministic unitary M\"obius function and discussed corresponding conjectures for the ordinary
M\"obius function. The deterministic main term and error for the squarefree pair count follow from the additive theorem for power-free integers of
Brandes \cite{Brandes}. See Lemma~\ref{lem:sf} for details.
To see how the distribution of signs on primes affects
densities and pair counts on the integers, first consider independent signs at the integers. If the signs are independent and
have the common distribution $\PP(f(n)=\sigma)=\pi_\sigma$, then for a fixed $N$ distinct unordered pairs $\{a,N-a\}$
use disjoint sets of random variables. Grouping the two ordered contributions from each unordered pair and applying Hoeffding's inequality
\cite[Theorem~2]{Hoeffding} and the first Borel--Cantelli lemma gives, almost surely for all sufficiently large integers,
\[
R_{\sigma\tau}(N)=\pi_\sigma\pi_\tau N+
O_{\mathcal X}\bigl(\sqrt{N\log N}\bigr).
\]
The central position when $N$ is even changes a count by at most one, which is absorbed into the error. For fair signs the main term is $N/4$.

This calculation gives a benchmark. With independent signs at the integers, the error is at most of order $\sqrt{N\log N}$.
In the multiplicative models the randomness comes from the primes. Even when two pairs involve four distinct integers,
those integers may have common prime factors. This raises two questions.
How do biases at the primes affect the main term and the error in the pair counts?
When the basic signs are fair or nearly fair, does multiplicative dependence still permit errors of order $\sqrt{N\log N}$?

\paragraph{Contributions}
\begin{itemize}
\item We study the four ordered sign counts in the reflection pairs $(f(a),f(N-a))$. The main term is $N/4$ in the nonvanishing model and $F(N)/4$ in the squarefree-supported model, where $F(N)$ counts all pairs with both summands squarefree.
\item Theorem~\ref{thm:main}\textup{(i)} estimates the errors for opposite-sign and same-sign pairs using the local weight $A_N$ and the full weight $B(N)$. If $A_N\to\infty$, each count in the nonvanishing model is $N/4+o(N)$. This includes constant negative-sign probability $0<\rho<1$.
\item For square-summable biases of the independent basic signs, Theorem~\ref{thm:main}\textup{(ii)} gives an error bound of order $N^{1/2+\varepsilon}$ in each of the four models in Table~\ref{tab:models}.
\item Under the same bias condition, Theorem~\ref{thm:cm-extrema} determines the optimal error scale $\sqrt{N\log N}$ and the exact constants in the completely multiplicative model. Each normalized counting error has limit superior $1/2$ and limit inferior $-1/2$. For the normalized reflection sum these limits are $2$ and $-2$. The extremes in both directions are approached in the initial interval $[X,X+X/\log X)$ of each dyadic block.
\end{itemize}

Balance of the marginal densities, even together with
$\sum_qm_q/q=\infty$, is insufficient to rule out reflection obstructions arising from long intervals of agreement with quadratic characters.
The random squarefree model provides a probabilistic comparison for M\"obius-type problems and a benchmark for their errors. Determining the actual values requires more precise calculations.
Section~\ref{sec:main} presents the main results and preliminary lemmas.
Section~\ref{sec:errors} proves the two types of error bounds.
Section~\ref{sec:cm-proof} determines the exact extremes in the completely multiplicative model.
Section~\ref{sec:scope} discusses a threshold for sparse biases,
pair distributions and density obstructions, the scope of the models, and open problems.
The fixed-bias estimates and the block counterexample are proved in Appendices~\ref{app:biased} and~\ref{sec:blocks}, respectively.

\subsection{Notation and probability-space conventions}
The letter $q$ always denotes a prime. We write $v_q(n)$ for the prime valuation, $\omega(n)$ for the number of distinct prime factors,
$\operatorname{rad}(n)=\prod_{q\mid n}q$ for the radical, $\tau(n)$ for the number of positive divisors,
and $\mu^2(n)$ for the squarefree indicator. All logarithms are natural.
On the probability space $(\Omega,\mathcal F,\PP)$ of the model at hand, $\mathcal X\in\Omega$
denotes one infinite realization of the random signs. For each fixed model and law, the countably many probability-$1$ events needed for the conclusions can be intersected. Denote the resulting event by $\Omega_0$. The notation $A(N)\ll_{\varepsilon,\mathcal X}B(N)$ or
$A(N)=O_{\varepsilon,\mathcal X}(B(N))$ means that there are a finite constant
$C(\varepsilon,\mathcal X)$ and a threshold $N_0(\varepsilon,\mathcal X)$ for which the bound holds
for all $N\ge N_0$. The constant may also depend on any other fixed parameters indicated.
We write $S_f(x)=\sum_{1\le n\le x}f(n)$ for the one-point partial sum.
In the nonvanishing model, $S_f(N)=S_f(N-1)+O(1)$.

\section{Main results}\label{sec:main}
\subsection{Statements of the results}
\begin{definition}[Prime-vector model]\label{def:model}
For each prime $q$, let $X_q=(X_{q,k})_{k\ge1}$ be a random vector with coordinates in $\pm1$.
The vectors at different primes are independent, and the coordinates at a given prime may have arbitrary dependence. Set
\[
\PP(X_{q,1}=-1)=\rho_q,\qquad
f(1)=1,\qquad f(n)=\prod_{q^k\parallel n}X_{q,k}.
\]
We call this the nonvanishing model, since $f(n)$ takes only the values $-1$ and $1$. Its squarefree-supported version is
\[
g(n)=\mu^2(n)\prod_{q\mid n}X_{q,1}.
\]
\end{definition}
Multiplicativity in the definition is required on coprime arguments. The completely multiplicative model is obtained by taking $X_{q,k}=X_{q,1}^k$. The basic sign at a prime is denoted by $\xi_q=X_{q,1}$. In the independent prime-power model we write $\xi_{q,k}=X_{q,k}$.
This is the standard definition of a multiplicative function by its values on prime powers \cite{DLMF}.

Set
\begin{equation}\label{eq:probability-scale}
\begin{split}
m_q&=\min(\rho_q,1-\rho_q),\qquad y_N=(\log N)^3,\\
A_N&=\sum_{y_N<q\le N}\frac{m_q}{q},\qquad B(x)=\sum_{q\le x}\frac{m_q}{q},\\
F(N)&=\sum_{a=1}^{N-1}\mu^2(a)\mu^2(N-a).
\end{split}
\end{equation}
Here $0\le m_q\le1/2$, and the mean of the first-power coordinate satisfies
$|\EE X_{q,1}|=|1-2\rho_q|=1-2m_q$.
All asymptotic statements concern sufficiently large $N$.

To use common notation for the two main terms and the errors of the four counts, set
\[
M_h(N)=\begin{cases}F(N),&h=g,\\N,&h\text{ is nonvanishing},\end{cases}
\quad
\Delta_h(N)=\max_{\sigma,\tau\in\{-1,1\}}
\left|R^h_{\sigma\tau}(N)-\frac{M_h(N)}4\right|.
\]
\begin{identitylemma}
In the nonvanishing model, set $C_f(N)=\sum_{a=1}^{N-1}f(a)f(N-a)$. Then
\begin{equation}\label{eq:identity}
4R^f_{\sigma\tau}(N)=N-1+(\sigma+\tau)S_f(N-1)+\sigma\tau C_f(N).
\end{equation}
In the squarefree-supported model, set
\[
U(N)=\sum_{a=1}^{N-1}g(a)\mu^2(N-a),\qquad
V(N)=\sum_{a=1}^{N-1}g(a)g(N-a),
\]
Then
\begin{equation}\label{eq:identity-sf}
4R^g_{\sigma\tau}(N)=F(N)+(\sigma+\tau)U(N)+\sigma\tau V(N).
\end{equation}
\end{identitylemma}

Expand $\ind_{\{f(n)=\sigma\}}=(1+\sigma f(n))/2$ and
$\ind_{\{g(n)=\sigma\}}=(\mu^2(n)+\sigma g(n))/2$, and substitute $a\mapsto N-a$ in the cross term.

The general error bound allows arbitrary dependence within each prime vector, as in Definition~\ref{def:model}. The near-fair error bound applies to the four rules in Table~\ref{tab:models}.
In the table, $\xi_j$ are the independent basic signs, and $j$ ranges over all the listed coordinates.

\begin{table}[H]
\centering
\caption{The four rules with independent basic signs}\label{tab:models}
\renewcommand{\arraystretch}{1.25}
\begin{tabular}{lll}
\hline
Model & Function or prime-power rule & Basic coordinate $j$\\
\hline
Squarefree-supported & $g(n)=\mu^2(n)\prod_{q\mid n}\xi_q$ & Primes $q$\\
Completely multiplicative & $f(q^k)=\xi_q^{\,k}$ & Primes $q$\\
Common sign on prime powers & $f(q^k)=\xi_q$ & Primes $q$\\
Independent prime powers & $f(q^k)=\xi_{q,k}$ & All $(q,k)$ with $k\ge1$\\
\hline
\end{tabular}
\end{table}

The first row is the squarefree-supported model. The last three rows are special cases of the nonvanishing model. In the common-sign model, all positive powers of a prime have the same sign.

The model is called fair if each basic sign is positive or negative with probability $1/2$ each.
It is called near-fair if the biases $\PP(\xi_j=-1)-1/2$ of the basic signs are square summable.

\begin{theorem}[General and near-fair error bounds]\label{thm:main}
\begin{enumerate}
\item[\textup{(i)}] Take the independent prime vectors in Definition~\ref{def:model}, allowing arbitrary dependence within each vector. On $\Omega_0$, for $h=f,g$, $\sigma=\pm1$ and every sufficiently large integer $N$,
\begin{equation}\label{eq:weight-mixed}
 \left|R^h_{\sigma,-\sigma}(N)-M_h(N)/4\right|\ll_{\mathcal X} Ne^{-4A_N},
\end{equation}
\begin{equation}\label{eq:weight-all}
 \Delta_h(N)
 \ll_{\mathcal X} N\bigl(e^{-4A_N}+e^{-2B(N)}\bigr).
\end{equation}
If $A_N\to\infty$, then $B(N)\ge A_N$ gives that both errors are $o(N)$.
\item[\textup{(ii)}] In any of the four models $h$ in Table~\ref{tab:models}, suppose that the independent basic signs satisfy
\[
\sum_j\left(\PP(\xi_j=-1)-\frac12\right)^2<\infty,
\]
Then on $\Omega_0$, for every $\varepsilon>0$ and all sufficiently large integers $N$,
\[
\Delta_h(N)\ll_{\varepsilon,\mathcal X}N^{1/2+\varepsilon}.
\]
In particular, the sum of squares is $0$ in the fair case, so the bound applies.
\end{enumerate}
\end{theorem}

Why treat opposite-sign and same-sign pairs separately? In the basic identities \eqref{eq:identity} and \eqref{eq:identity-sf},
the coefficient of the one-point term is $\sigma+\tau$. Thus the reflection sum alone controls the opposite-sign counts, while the same-sign counts involve both terms.
The one-point bound uses the full weight $B(N)$. If $B(N)\ge2A_N$,
its contribution to the bound in \eqref{eq:weight-all} is no greater than that of the reflection sum.

\begin{theorem}[Optimal error in the completely multiplicative model]\label{thm:cm-extrema}
Let $f$ be completely multiplicative. Suppose that the prime signs $\xi_q=f(q)$ are independent and that
$\sum_q(\PP(\xi_q=-1)-1/2)^2<\infty$.
\[ s_N=\sqrt{N\log N},\qquad D_{\sigma\tau}(N)=\frac{R^f_{\sigma\tau}(N)-N/4}{s_N}.
\]
On $\Omega_0$, for all $\sigma,\tau\in\{-1,1\}$, the errors of the four counts satisfy
\[
 \limsup_{N\to\infty}D_{\sigma\tau}(N)=\tfrac12,\qquad
 \liminf_{N\to\infty}D_{\sigma\tau}(N)=-\tfrac12,\qquad
 \limsup_{N\to\infty}\frac{\Delta_f(N)}{s_N}=\tfrac12 .
\]
The corresponding reflection sum satisfies
\begin{equation}\label{eq:cm-explicit}
 \limsup_{N\to\infty}\frac{C_f(N)}{s_N}=2,\qquad
 \liminf_{N\to\infty}\frac{C_f(N)}{s_N}=-2 .
\end{equation}
Along $X=2^j$, let $I_X=[X,X+X/\log X)\cap\mathbb Z$. On $\Omega_0$, normalize the quantities for each $N$ in $I_X$ before taking their maxima and minima. The resulting limits are
\begin{equation}\label{eq:cm-four-counts}
 \begin{aligned}
 \lim_{j\to\infty}\left(\max_{N\in I_X}D_{\sigma\tau}(N),\min_{N\in I_X}D_{\sigma\tau}(N)\right)
   &=(1/2,-1/2),\\
 \lim_{j\to\infty}\left(\max_{N\in I_X}\frac{C_f(N)}{s_N},\min_{N\in I_X}\frac{C_f(N)}{s_N}\right)
   &=(2,-2).
 \end{aligned}
\end{equation}
\end{theorem}

On $\Omega_0$, for every $\eta>0$ and all sufficiently large $N$, we have
$\Delta_f(N)\le(1/2+\eta)s_N$. For any $0<\eta<1/2$, in each sufficiently large initial dyadic interval $I_X$,
each count has some $N_+$ for which $D_{\sigma\tau}(N_+)\ge1/2-\eta$
and some $N_-$ for which $D_{\sigma\tau}(N_-)\le-(1/2-\eta)$. The witnessing integers may depend on the sign pair and on the direction of the deviation.
Thus $\sqrt{N\log N}$ is the optimal error scale, and $1/2$ is the best asymptotic constant for the absolute errors of the four counts.
The limit superior and limit inferior along the full sequence determine the constants. The short-interval extremes show that deviations in both directions occur near the start of every sufficiently large dyadic block.

The two theorems use progressively more structure. Theorem~\ref{thm:main}\textup{(i)} allows arbitrary dependence among the powers of a given prime.
Part~\textup{(ii)} uses Walsh moment estimates for independent basic signs in four specified models to obtain the square-root-scale bound.
The completely multiplicative model has additional structure. Each term $f(a)f(N-a)$ depends only on the squarefree kernel of $a(N-a)$.
This controls correlations between different targets and lets us use relative moderate deviations to show that the reflection sum reaches the scale of the upper bound.
Atherfold's one-point bound \cite[Theorem~1]{AtherfoldExtended} gives $S_f(N-1)=o(s_N)$, so the reflection sum determines the extremes of the four counts.

\begin{remark}[Scope of the concentration estimates]
In a multiplicative model, changing a basic coordinate at the prime $q$ can affect $O(N/q)$ pair positions.
Applying bounded differences directly to all primes gives an upper bound of $O(N^2)$ for the sum of squared differences.
For the square-root-scale error in Theorem~\ref{thm:main}\textup{(ii)}, see the Walsh moment estimates in Section~\ref{sec:fair-proof}.
The two-sided extremes in the completely multiplicative model are proved in Section~\ref{sec:cm-proof}.
\end{remark}

Lemma~\ref{lem:sf} gives the positive linear lower bound $F(N)\ge((10-\pi^2)/4)N-1$ for the squarefree pair count.
Hence if $A_N\to\infty$ in Theorem~\ref{thm:main}\textup{(i)}, or under the assumptions of part~\textup{(ii)},
almost surely each of the four sign patterns represents every sufficiently large integer in the corresponding model.
\begin{remark}[Deterministic error in the squarefree main term]
With the notation for the singular series in Lemma~\ref{lem:sf}, Theorem~\ref{thm:main}\textup{(i)} gives
\[
R^g_{\sigma\tau}(N)=\frac{\mathfrak S_{\mathrm{sf}}(N)}4N
+O_{\varepsilon,\mathcal X}\!\left(Ne^{-4A_N}+
 \ind_{\{\sigma=\tau\}}Ne^{-2B(N)}+N^{7/11+\varepsilon}\right).
\]
For the squarefree model in Theorem~\ref{thm:main}\textup{(ii)}, the corresponding estimate is
\[
R^g_{\sigma\tau}(N)=\frac{\mathfrak S_{\mathrm{sf}}(N)}4N
+O_{\varepsilon,\mathcal X}(N^{7/11+\varepsilon}).
\]
Here $0<\varepsilon<4/11$ can be any fixed value in the indicated range. The square-root-scale bound about the exact center $F(N)/4$ still holds.
Replacing that center by the Euler-product main term also introduces the deterministic error in the squarefree pair count.
\end{remark}

\begin{corollary}[Fixed bias and fair signs]\label{cor:fixed}
Fix $0<\rho<1$. If $\rho_q=\rho$ and $m=\min(\rho,1-\rho)$, the two models in Definition~\ref{def:model} almost surely satisfy
\[
\Delta_f(N)+\Delta_g(N)\ll_{\rho,\mathcal X}N/(\log N)^{2m}.
\]
In the completely multiplicative model, the maximum error has the following bounds in the different ranges.
\begin{enumerate}
\item If $0<\rho<2/3$ and $\rho\ne1/2$, then
\[
\Delta_f(N)\sim\frac{|c_{\rho,f}|}{2}\frac{N}{(\log N)^{2\rho}}.
\]
The nonzero random constant $c_{\rho,f}\ne0$ is specified in Appendix~\ref{app:biased}.
\item If $\rho=1/2$, then simultaneously for every $\varepsilon>0$,
$\Delta_f(N)\ll_{\varepsilon,f}N^{1/2+\varepsilon}$.
This bound also holds for the remaining models in Table~\ref{tab:models} when their basic signs are fair. The exact extremes for the fair completely multiplicative model are given in Theorem~\ref{thm:cm-extrema}.
\item If $2/3\le\rho<1$, then
\[
\frac{N}{(\log N)^{2\rho}}\ll_{\rho,f}\Delta_f(N)
\ll_{\rho,f}N\left(\frac{\log\log N}{\log N}\right)^{4m}.
\]
\end{enumerate}
In addition, for each fixed $0<\rho<1$, the completely multiplicative reflection sum almost surely satisfies
\[
|C_f(N)|:=\left|\sum_{a=1}^{N-1}f(a)f(N-a)\right|
\ll_{\rho,f}N\left(\frac{\log\log N}{\log N}\right)^{4m}.
\]
For $2/3\le\rho<1$, there is still a gap between the upper and lower bounds.
\end{corollary}
\begin{proof}
By Mertens' second theorem \eqref{eq:mertens-reciprocal},
\[
B(N)=m\log\log N+O(1),\qquad A_N=B(N)-m\log\log y_N+O(1).
\]
Substituting these estimates into Theorem~\ref{thm:main}\textup{(i)} gives the general bound. The reflection-sum estimate is given in \eqref{eq:weight-reflection}.
Appendix~\ref{app:biased} gives the asymptotic for the one-point sum in the completely multiplicative model and its contribution to the maximum error.
\end{proof}
\begin{remark}[Mechanisms for fair signs and fixed bias]
In the completely multiplicative model for $\rho\ne1/2$, the Selberg--Delange method gives the main term
$c_{\rho,f}x/(\log x)^{2\rho}$ for the one-point sum, where
$c_{\rho,f}=H_f/\Gamma(1-2\rho)$ and $H_f>0$. See Appendix~\ref{app:biased}.
The logarithmic power corresponds to the dominant factor of type $\zeta(s)^{1-2\rho}$ in the Dirichlet series near $s=1$.
At the fair point the reciprocal gamma factor vanishes at $0$. That cancellation alone does not imply a square-root-scale error.
See Section~\ref{sec:fair-proof} for the fair-model error bounds and Section~\ref{sec:cm-proof} for the exact extremes in the completely multiplicative model.
The almost-sure conclusions of the corollary hold for each fixed $\rho$. If any prescribed countable family of parameters is realized on the same probability space
(for example, rational $\rho$), one can intersect the corresponding probability-$1$ events. A uniform result over uncountably many parameters would require a new argument.
\end{remark}

\begin{samepage}
\begin{theorem}[Threshold and error for sparse biases]\label{thm:sparse}
Fix $\alpha>0$. Take independent prime signs and set
\begin{equation}\label{eq:family}
\rho_q=\min\left\{\frac12,
 [\log\log(q+e^e)]^{-\alpha}\right\},\qquad X_{q,k}=1\quad(k\ge2).
\end{equation}
Then the relative density of negative primes is almost surely $0$. For $0<\alpha\le1$, define
\[
E_\alpha(N)=\begin{cases}
N\exp\!\left(-\frac{2(\log\log N)^{1-\alpha}}{1-\alpha}\right),&0<\alpha<1,\\
N/(\log\log N)^2,&\alpha=1.
\end{cases}
\]
Almost surely the following uniform error bound holds.
\[
\Delta_f(N)+\Delta_g(N)\ll_{\alpha,\mathcal X}E_\alpha(N).
\]
If $\alpha>1$, then almost surely $\Delta_f(N)=o(N)$ fails.
\end{theorem}
\end{samepage}
\noindent\textbf{Models and error scales.}
The table summarizes bounds in common cases. The counts are centered at $M_h(N)/4$. Unless otherwise indicated, the errors hold on the same realization for every sufficiently large integer target.
\begin{center}
\small
\begin{tabular}{@{}>{\raggedright\arraybackslash}p{.36\linewidth}l>{\raggedright\arraybackslash}p{.43\linewidth}@{}}
\hline
Model and conditions & Result & Error scale\\
\hline
Independent prime vectors, including the squarefree-supported version &
Thm.~\ref{thm:main}(i) &
Opposite signs $Ne^{-4A_N}$, same signs $N(e^{-4A_N}+e^{-2B(N)})$\\
Four models, square-summable biases of the basic signs &
Thm.~\ref{thm:main}(ii) & $N^{1/2+\varepsilon}$\\
Near-fair completely multiplicative model &
Thm.~\ref{thm:cm-extrema} & $\sqrt{N\log N}$, optimal asymptotic constant $1/2$\\
General model, fixed $0<\rho<1$ &
Cor.~\ref{cor:fixed} & $N/(\log N)^{2\min(\rho,1-\rho)}$\\
Sparse biases in \eqref{eq:family}, $0<\alpha\le1$ &
Thm.~\ref{thm:sparse} & $E_\alpha(N)$\\
\hline
\end{tabular}
\end{center}

These results show that the distribution of negative-sign probabilities at the primes affects the error bounds. It remains to identify how that distribution determines the error scale and to find a general threshold for approximate equidistribution of the four sign patterns.
\begin{theorem}\label{thm:counter}
There exist probabilities $\rho_q\in[0,1]$ such that
\[
\sum_q\frac{m_q}{q}=\infty,
\]
Write $H$ for the completely multiplicative function generated by these probabilities, and $f$ for the non-completely multiplicative function constructed with the same probabilities. They have the following properties.
\begin{enumerate}
\item[\textup{(i)}] \textbf{Completely multiplicative counterexample.} For every realization of $H$, there are infinitely many primes $P\equiv7\pmod8$ for which
\[
R^H_{--}(P)=R^H_{++}(P)=0.
\]
\item[\textup{(ii)}] \textbf{Non-completely multiplicative counterexample.} The function $f$ is a genuinely non-completely multiplicative, nonvanishing model. At these targets,
\[
C_f(P):=\sum_{a=1}^{P-1}f(a)f(P-a)=-(P-3)/2,
\]
so the quarter law fails.
\item[\textup{(iii)}] \textbf{Additional properties.} Almost surely, both $f$ and $H$ have mean zero. For every fixed real Dirichlet character $\chi$, we have $\sum_q(1-f(q)\chi(q))/q=\infty$, and the same holds for $H$.
For every $0<\kappa<1$ one can also choose a family of probabilities with $0<\rho_q<1$ such that both failures in (i) and (ii) occur infinitely often on a common event of probability at least $1-\kappa$.
\end{enumerate}
\end{theorem}
\begin{remark}
Theorem~\ref{thm:counter} shows that divergence of the total probability weight $\sum_qm_q/q=\infty$ alone does not guarantee the quarter law.
The long character-matching intervals in the construction give $0\le A_P\le1/(2P)$ at certain targets. This illustrates the role of the scale assumption.
The construction does not determine whether $A_N\to\infty$ is necessary for the quarter law.
The completely multiplicative version rules out same-sign representations at the specified targets. The non-completely multiplicative version gives asymptotic imbalance.
\end{remark}

\subsection{Preliminary lemmas}

We use the following three Mertens formulas for primes. See Apostol \cite[Sections~4.6--4.8]{Apostol}.
\begin{align}
 &\sum_{p\le x}\frac{\log p}{p}=\log x+O(1)
   &&\text{(Mertens' first theorem)},\label{eq:mertens-weighted}\\
 &\sum_{p\le x}\frac1p=\log\log x+O(1)
   &&\text{(Mertens' second theorem)},\label{eq:mertens-reciprocal}\\
 &\prod_{p\le x}\left(1-\frac1p\right)^{-1}\asymp\log x
   &&\text{(order of magnitude in Mertens' third theorem)}.\label{eq:mertens-product}
\end{align}
The third formula also follows by taking logarithms of the product and using the second formula and $\sum_p p^{-2}<\infty$.
\begin{lemma}\label{lem:weight-mean}
For every multiplicative function with $0\le v(n)\le1$ and $v(1)=1$, the following bound holds uniformly.
\[
 \sum_{n\le x}v(n)\ll\sqrt x+x\exp\!\left(-\sum_{q\le x}\frac{1-v(q)}q\right)
 \qquad(x\ge3).
\]
\end{lemma}
\begin{proof}
Substitute $\log n=\sum_{q^k\parallel n}k\log q$, group the terms by $n=aq^k$, and drop the restriction $(a,q)=1$ to obtain
\[
 \sum_{n\le x}v(n)\log n\le\sum_{a\le x}v(a)\sum_{q^k\le x/a}k\log q
 \ll x\sum_{a\le x}v(a)/a .
\]
Chebyshev's upper bound $\sum_{q\le t}\log q\ll t$ controls the first powers of primes. The higher powers contribute
$O(\sqrt t\log^3(2t))=O(t)$. Split the sum at $\sqrt x$ and bound the harmonic sum by an Euler product.
Each factor is at most $1+v(q)/q+O(1/q^2)$, and Mertens' second theorem \eqref{eq:mertens-reciprocal} gives the bound.
\end{proof}

\begin{lemma}\label{lem:McD}
Put $\mathcal F_y=\sigma(X_q:q\le y)$. Conditional on $\mathcal F_y$, suppose that a real random variable $Z$ depends only on finitely many of the remaining prime vectors.
Changing the vector at the prime $q$ changes $Z$ by at most $c_q$. If $\sum c_q^2>0$, then
\[
 \PP\bigl(|Z-\EE(Z\mid\mathcal F_y)|>t\mid\mathcal F_y\bigr)
 \le2\exp\!\left(-2t^2/\sum_{q>y}c_q^2\right)\qquad(t>0).
\]
\end{lemma}
\begin{proof}
Apply McDiarmid's inequality \cite{McDiarmid}.
\end{proof}

\begin{lemma}[Lower bound and main term for squarefree pairs]\label{lem:sf}
For $N\ge2$, we have
\[
F(N)\ge c_*N-1,\qquad c_*=(10-\pi^2)/4>0.
\]
In addition, put
\begin{equation}\label{eq:sf-series}
\begin{split}
\mathfrak S_{\mathrm{sf}}(N)
&=\prod_{q^2\mid N}\left(1-\frac1{q^2}\right)
  \prod_{q^2\nmid N}\left(1-\frac2{q^2}\right)\\
&=C_{\mathrm{sf}}\prod_{q^2\mid N}\frac{q^2-1}{q^2-2},
\qquad C_{\mathrm{sf}}=\prod_q\left(1-\frac2{q^2}\right).
\end{split}
\end{equation}
Then for every fixed $\varepsilon>0$ the following unconditional estimate holds uniformly in integer $N$.
\begin{equation}\label{eq:sf-asymptotic}
F(N)=\mathfrak S_{\mathrm{sf}}(N)N+O_\varepsilon(N^{7/11+\varepsilon}).
\end{equation}
\end{lemma}
\begin{proof}
If either summand is not squarefree, some prime square divides it. A union bound over the two summands gives
\[
F(N)\ge N-1-2N\sum_q q^{-2}.
\]
Also,
\[
\sum_q q^{-2}\le\frac14+\sum_{j\ge1}(2j+1)^{-2}
=\frac{\pi^2}{8}-\frac34,
\]
Substitution gives the lower bound for each $N$.

For the main term, apply Brandes' theorem on sums and differences of power-free numbers \cite[Theorem~1]{Brandes} with
$k=l=2$ and $x=n=N$. The additive count in that theorem requires $a\le x$ and both $a$ and $n-a$ to be positive and power-free with their respective exponents.
Here $a\le x$ is automatic, and positivity of both summands gives exactly $1\le a<N$, so $R^+_{2,2}(N;N)=F(N)$.
The main term uses $X^+(x;n)=\min(x,n)$, which here equals $N$, and
$C_{2,2}(N)S_{2,2}=\mathfrak S_{\mathrm{sf}}(N)$.
Since $\tfrac43(k-1)\le l\le4k$, the error exponent in that theorem is
\[
\frac{7(k+l)}{6k+2l+7kl}\bigg|_{k=l=2}
=\frac{28}{44}=\frac7{11}.
\]
In the error, $x+n=2N$ is of the same order as $N$. The implied constant in the theorem is independent of $x,n$,
so \eqref{eq:sf-asymptotic} holds uniformly over all integer $N$.
The local factors in the main term come from the forbidden classes $a\equiv0,N$ modulo $q^2$.
The classes coincide when $q^2\mid N$ and are distinct otherwise.
This gives $1-q^{-2}$ and $1-2q^{-2}$, respectively, as in \eqref{eq:sf-series}.
The product converges absolutely and is strictly positive because the deviations of the local factors from $1$ are absolutely summable.
The local density is determined by the divisibility condition $q^2\mid N$. The error is controlled by Brandes' theorem.
\end{proof}

This is the same count as the quantity in Tian and Yang \cite[Definition 4.1 and Proposition 4.2]{TianYang}
denoted by $M_2(N)$. Their paper gives an elementary error bound $O_\varepsilon(N^{2/3+\varepsilon})$.
Here we use Brandes' stronger result.

\section{Proof of the general error bounds}\label{sec:errors}
\paragraph{Proof outline}
For the probability-weight bounds, we use bounded differences to control fluctuations, an upper-bound sieve for the conditional mean of the reflection sum, and a sequence of decreasing cutoffs for the one-point sum.
For the near-fair bounds, we check the Walsh coefficients and moment estimates under the fair law in each of the four models, then transfer the almost-sure conclusions by absolute continuity.
\subsection{Probability-weight error bounds}\label{sec:scale}
We prove Theorem~\ref{thm:main}\textup{(i)}. The counting identities reduce the problem to two types of sums.
Both factors in a reflection term have random signs, while only one factor in a one-point term does. We first fix the values at small primes
and use bounded differences to control fluctuations about the conditional mean. We then estimate the two types of conditional means separately. Write
\[
 P_0=\lfloor(\log(2N))^5\rfloor,\qquad
 W(N,P)=\exp\bigl(-2(B(N)-B(P))\bigr).
\]
By $2m_q\le1$ and Mertens' second theorem \eqref{eq:mertens-reciprocal}, uniformly for $2\le P\le N$ we have
$W(N,P)\gg\log P/\log N$ and $e^{-2B(N)}\gg1/\log N$.

\subsubsection{Conditional fluctuations and the reflection sum}
First fix the prime vectors at primes up to $P_0$. Write
\[
 Z_{\rm ref}\in\{C_f(N),V(N)\},\qquad
 Z_{\rm one}\in\{S_f(N-1),U(N)\}.
\]
Changing the vector at $q>P_0$ affects at most $2N/q$ positions in the reflection sum and changes that sum by at most $4N/q$.
In the one-point sum at most $N/q$ random factors are affected, and the change is at most $2N/q$.
The sum of squared differences is $\ll N^2/P_0$ for each type. By Lemma~\ref{lem:McD}, for a sufficiently large absolute constant $K$,
the conditional failure probability at the threshold $KN\sqrt{\log N/P_0}$ is $O(N^{-3})$ in each case. Removing the conditioning and summing over the four sums and all $N$,
Borel--Cantelli gives the following eventual estimate on $\Omega_0$.
\begin{equation}\label{eq:weight-outer}
 |Z_{\rm ref}-\EE(Z_{\rm ref}\mid\mathcal F_{P_0})|
 + |Z_{\rm one}-\EE(Z_{\rm one}\mid\mathcal F_{P_0})|
 \ll N/\log^2 N .
\end{equation}

To estimate the conditional mean of $Z_{\rm ref}$, put $z=N^{1/4}$. If $P_0<q\le z$ and $q\nmid N$,
the absolute local mean is $1-2m_q$ at positions where $q\parallel a$ or $q\parallel N-a$.
The two summands correspond to the distinct residue classes $0,N$ modulo $q$.
Each class has $q$ lifts modulo $q^2$. Exactly one of those lifts makes the corresponding summand divisible by $q^2$.
Removing these two lifts leaves $2(q-1)$ classes. At all other positions, we bound the local factor by $1$.
If $q\mid N$, the two zero classes coincide modulo $q$, so we bound the local factor for that prime by $1$ everywhere.
In both cases it suffices that the absolute value of any product of prime-power signs is at most $1$.
Independence across primes therefore bounds the absolute conditional mean of the reflection sum
by the sum of these products of local weights over $1\le a<N$.

To estimate this sum, use an auxiliary space independent of the original model. For each $P_0<q\le z$ with $q\nmid N$,
select the prime independently with probability $2m_q$. If selected, forbid the $2(q-1)$ residue classes modulo $q^2$ described above.
For a fixed $a$, the survival probability is exactly the product of its local weights. Fix one outcome of the selection. Write
$\mathcal A=\{1,\ldots,N-1\}$ and let $\mathcal A_q$ consist of the positions in the forbidden classes for $q$.
For unselected primes and primes dividing $N$, set $\mathcal A_q=\varnothing$.
Let $\gamma(q)$ be the density of the forbidden classes modulo $q^2$, and extend it multiplicatively to squarefree $d$.
The Chinese remainder theorem gives
\[
 \#\bigcap_{q\mid d}\mathcal A_q=(N-1)\gamma(d)+r_d,\qquad
 |r_d|\le d\,2^{\omega(d)} ,
\]
because the counting error in each class modulo $d^2$ is at most $1$.
Also $0\le\gamma(q)\le\min(2/q,1/2)$, so the sieve-density condition in Ford
\cite[Theorem~2.4 and Exercise~2.1]{FordSieve} holds in the form
\[
 \prod_{u\le q\le v}(1-\gamma(q))^{-1}
 \le (\log v/\log u)^2\exp(O(1/\log u))\qquad(2\le u\le v\le z),
\]
with constants independent of the selected classes and of $N$. The theorem bounds the number of survivors by
\[
 \ll (N-1)\prod_{q\le z}(1-\gamma(q))
   +\sum_{d\le z}\mu^2(d)|r_d|.
\]
Since $2^{\omega(d)}\le\tau(d)$ and $\sum_{d\le z}\tau(d)\ll z\log(2z)$,
the sum of the remainders is $\ll z^2\log(2z)$. Averaging over the auxiliary selection gives the survival-density product
\[
 \prod_{\substack{P_0<q\le z\\q\nmid N}}
       \left(1-\frac{4m_q(q-1)}{q^2}\right)\ll e^{-4A_N}.
\]
The sums of reciprocal primes from $y_N$ to $P_0$ and from $z$ to $N$ are bounded.
Removing positions of square divisibility from each set of forbidden classes costs only
$O(\sum_q m_q/q^2)=O(1)$ in the exponent. Omitting primes dividing $N$ costs
$O(\sum_{q\mid N,\,q>P_0}m_q/q)\ll\log N/P_0$. Combining the sieve remainder with \eqref{eq:weight-outer} and using
$e^{-4A_N}\gg1/\log^2 N$ gives
\begin{equation}\label{eq:weight-reflection}
 |C_f(N)|+|V(N)|\ll_{\mathcal X}Ne^{-4A_N}.
\end{equation}

\subsubsection{One-point sums and the full prime weight}
Let $h=f$ or $g$ and write $h_P(a)=\EE(h(a)\mid\mathcal F_P)$. In the prime decomposition we use
the actual signs at $p\le P$, retaining the squarefree restriction for $g$, and the mean of each relevant prime power at the remaining primes. Set
\[
 H_P^f(N)=\sum_{a<N}f_P(a),\qquad H_P^g(N)=\sum_{a<N}g_P(a)\mu^2(N-a),
\]
Here and below the summation positions are positive integers. In $H_P^g(N)$, the factor $g_P(a)$ retains the squarefree restriction on the first summand, and the deterministic mask $\mu^2(N-a)$ restricts the second.

For an integer cutoff $P$ set $L=\log P,\ b=\lfloor L^{30}\rfloor,\ D=L^2,\ T=e^{L^2}$.
Take $P$ large enough that $2\le b\le\sqrt P$ and $D^2<P$. We call positive integers whose prime factors are at most $P$
$P$-smooth, including $1$. Write $h_b(a)=\EE(h(a)\mid\mathcal F_b)$, where $\mathcal F_b=\sigma(X_q:q\le b)$. On $\Omega_0$, for every sufficiently large integer $P$, every $1\le t\le2T$,
every modulus $k\le D^2$ and every residue class $r$, both models satisfy
\begin{equation}\label{eq:weight-smooth}
 \left|\sum_{\substack{d<t,\ d\equiv r\bmod k\\ d\text{ is }P\text{-smooth}}}
             (h(d)-h_b(d))\right|\ll tL^{-8}.
\end{equation}
Conditional on $\mathcal F_b$, the sum is centered. Changing one vector at $b<p\le P$ affects the sum by at most $2t/p$,
so the sum of squared differences is $\ll t^2/b$. Lemma~\ref{lem:McD} bounds each failure probability by $\le2\exp(-c b/L^{16})$.
A union bound over integer endpoints up to $\lceil2T\rceil$, with rounding up at most doubling the threshold, and the $O(L^8)$ residue classes
gives a total failure probability $\ll e^{L^2}L^8e^{-c'L^{14}}$, which is summable over all integers $P$. This proves \eqref{eq:weight-smooth}
on an event shared by all targets $N$.

On this event, when $P\le P_0$ and $P$ is sufficiently large, we have the one-step estimate
\begin{equation}\label{eq:weight-step}
 |H_P^h(N)-H_b^h(N)|\ll NW(N,P)L^{-2}.
\end{equation}
We first show that without the squarefree mask, the error in any one residue class modulo $k\le D^2$ is $O(NWL^{-7})$,
where $W=W(N,P)$. Write $a=md$ uniquely, with $d$ the full $P$-smooth part of $a$. Then
\[
 h_P(a)-h_b(a)=\EE h(m)\,[h(d)-h_b(d)].
\]
All prime factors of $m$ exceed $P$, while $k\le D^2<P$, so $\gcd(m,k)=1$.
If $a\equiv r\bmod k$ then $d\equiv m^{-1}r\bmod k$. Let $\Psi(t,P)$ count the integers up to $t$ that are
$P$-smooth. With $\alpha=1/L$, Mertens' first theorem \eqref{eq:mertens-weighted}
and his third theorem \eqref{eq:mertens-product} give
\[
 E_P=\prod_{p\le P}(1-p^{-1+\alpha})^{-1}\ll L,\qquad
 \Psi(t,P)\le t^{1-\alpha}E_P .
\]
The first bound uses $\sum_{p\le P}(p^\alpha-1)/p\ll1$. For the second, bound the contribution of each $d\le t$ by
$(t/d)^{1-\alpha}$. Hence the total contribution from $m<\sqrt N$
$\ll\sum_{m<\sqrt N}(N/m)^{1-\alpha}L\ll N^{1-\alpha/2}L\log N\ll NW L^{-7}$,
where $L\le5\log\log(2N)$.

For $m\ge\sqrt N$, group the terms by $J\le N/m<2J$, with $J$ ranging over nonnegative integer powers of $2$.
The outer sum in each group only involves $N/(2J)<m\le N/J$, all of whose prime factors exceed $P$.
Extending the range to $m\le N/J$ and retaining the restriction on prime factors gives
$\sum_m|\EE h(m)|\ll (N/J)W/L$. To see this, apply Lemma~\ref{lem:weight-mean} with weight zero
on positive powers of $p\le P$ and weight $|\EE h(p^\nu)|$ on the rest.
At the upper limit $N/J\ge\sqrt N$, the small primes give a factor $O(1/L)$, and $\sum_{\sqrt N<p\le N}1/p=O(1)$.
The square-root remainder is absorbed by $W/L\gg1/\log N$. In the $O(L^2)$ groups with $J\le T$,
the inner sum in \eqref{eq:weight-smooth} is at most $O(JL^{-8})$, for a total of $O(NWL^{-7})$.
For $J>T$ the inner sum is at most $2\Psi(2J,P)$, and
\[
 \sum_{J>T}\frac{\Psi(2J,P)}{J}
 \ll\sum_{d\text{ is }P\text{-smooth}}\min(T^{-1},d^{-1})
 \le T^{-\alpha}E_P\ll L/P .
\]
The first inequality follows by summing a dyadic geometric series for each $d$. The tail is therefore at most $O(NW/P)$, proving the estimate for each residue class.
Taking $k=1$ gives \eqref{eq:weight-step} for $f$.

For $g$, expand $\mu^2(N-a)=\sum_{\ell^2\mid N-a}\mu(\ell)$.
The terms with $\ell\le D$ give at most $D$ residue classes of moduli $\ell^2$, contributing $O(NWL^{-5})$ in total.
Set $E=(\log N)^2$. When $D<\ell\le E$, on the progression $a\equiv N\bmod\ell^2$,
both $|g_P(a)|,|g_b(a)|$ are bounded by $\prod_{p>P,\,p\mid a}(1-2m_p)$, since nonsquarefree positions contribute zero.
Retain only the factors with $P<p\le N^{1/4}$ and $p\nmid\ell$. In the index of the progression there is exactly one forbidden class for each $p$.
Select this class independently with probability $2m_p$. Given the selection, the intersection of the forbidden classes for squarefree $d$
contains at most one residue class modulo $d$ in the index of the progression. Thus $|r_d|\le1$, and the total remainder is $\le N^{1/4}$.
The local density is $0$ or $1/p$, so the same upper-bound sieve applies with dimension $1$. After averaging over the selection,
omitting primes with $p\mid\ell$ enlarges the exponential bound by a factor of at most $\exp(\sum_{p\mid\ell}2m_p/p)\le\ell/\varphi(\ell)$.
Thus the absolute contribution from each progression is
\[
 \ll (N/\ell^2+1)W\,\ell/\varphi(\ell)+N^{1/4}.
\]
Use the identity $n/\varphi(n)=\sum_{d\mid n}\mu^2(d)/\varphi(d)$. The convergence of the Euler product
$\sum_d\mu^2(d)/(d\varphi(d))<\infty$ gives $\sum_{\ell\le t}\ell/\varphi(\ell)\ll t$. 

Partial summation gives $\sum_{\ell>D}(\ell/\varphi(\ell))/\ell^2\ll1/D$. For $\ell>E$,
counting square divisors directly for $0<N-a<N$ shows that the contribution is $\ll N\sum_{\ell>E}\ell^{-2}\ll N/E$.
Thus the mask tail is at most $O(NW/D+E N^{1/4}+N/E)$, which is absorbed into $O(NW/L^2)$ by $W\gg L/\log N$.
This proves \eqref{eq:weight-step}.

Starting at $P_0$, take $P_{i+1}=\lfloor(\log P_i)^{30}\rfloor$ repeatedly until the first $P_s<P_*$.
Here $P_*$ is the finite, sample-dependent threshold above which the preceding estimates hold for all $P\ge P_*$. On this common event,
the error at each step, divided by $Ne^{-2B(N)}$, is at most
\[
 C\,e^{2B(P_i)}/\log^2P_i\ll1/\log P_i .
\]
Since $P_{i+1}\le\sqrt{P_i}$, the sum of these reciprocals is at most twice the last one. For the last level,
apply Lemma~\ref{lem:weight-mean} with weight $1$ on the powers of primes $p\le P_s$ and the absolute values of the corresponding means at the other prime powers,
giving $|H_{P_s}^h(N)|\ll\sqrt N+NW(N,P_s)\ll_{\mathcal X}Ne^{-2B(N)}$.
Sum over the levels and use \eqref{eq:weight-outer} to obtain
\begin{equation}\label{eq:weight-single}
 |S_f(N-1)|+|U(N)|\ll_{\mathcal X}Ne^{-2B(N)} .
\end{equation}

On $\Omega_0$, substitute \eqref{eq:weight-reflection} and \eqref{eq:weight-single}
into the basic identities \eqref{eq:identity} and \eqref{eq:identity-sf}. Only the reflection term remains for opposite signs.
For the same signs there is also $2\sigma S_f(N-1)$ or $2\sigma U(N)$. Changing the nonvanishing center from $(N-1)/4$ to $N/4$
costs only $1/4$. The squarefree center is exactly $F(N)/4$. This proves \eqref{eq:weight-mixed} and \eqref{eq:weight-all}.

\begin{remark}
Choose $y$ so that $\log N/y\to0$ and $\sum_{y<q\le N}m_q/q\to\infty$.
In the reflection proof, one may condition directly on primes up to $y$. The sum of reciprocals of primes shared by the summands is $\ll\log N/y$,
and the fluctuation threshold $KN\sqrt{\log N/y}=o(N)$ still gives a summable sequence of failure probabilities. Together with \eqref{eq:weight-single},
this gives the quarter-law asymptotic.
\end{remark}

\subsection{Near-fair error bounds}\label{sec:fair-proof}
We prove Theorem~\ref{thm:main}\textup{(ii)}. First take independent fair basic signs.
Distinct Walsh monomials are then orthogonal, so the second moment can be read directly from the coefficients.
We check the coefficients in each of the four models and use higher moments to obtain summable tail probabilities. Finally we handle bias by absolute continuity.

\subsubsection{Reduction to degree and moment bounds}
Write $(\xi_j)$ for the basic signs, indexed by primes or by prime powers.
For a finite coordinate set $L$, set $\xi_L=\prod_{j\in L}\xi_j$, with the empty set corresponding to the constant $1$.
For prime coordinates we also index monomials by squarefree integers $t$, writing $\xi_t=\prod_{q\mid t}\xi_q$.
These Walsh monomials are mutually orthogonal, so
\[
\left\|\sum_L c_L\xi_L\right\|_2^2=\sum_L|c_L|^2.
\]
If $n\le N$ has $r$ distinct prime factors, then $n\ge(r+1)!$. Hence
\[
\max_{n\le N}\omega(n)=O(\log N/\log\log N).
\]
For fair basic signs in the four models, $f(n)$ or a nonzero $g(n)$ is a monomial of degree at most $\omega(n)$.
Canceling repeated coordinates in a product by $\xi_j^2=1$ can only lower the degree.
Thus every statistic below is multilinear of degree at most
$d_N=O(\log N/\log\log N)$, and only finitely many variables occur.

\begin{lemma}\label{lem:walsh}
Let $Q_N$ be a finite multilinear polynomial in independent fair signs, of degree
$O(\log N/\log\log N)$. If for every $\delta>0$ it satisfies
\begin{equation}\label{eq:fair-L2}
\|Q_N\|_2^2\ll_\delta N^{1+\delta}.
\end{equation}
then on a single probability-$1$ event, for every $\varepsilon>0$,
$Q_N=O_{\varepsilon,\mathcal X}(N^{1/2+\varepsilon})$.
In particular, if $Q_N=\sum_Lm_L\xi_L$ with $m_L\ge0$ and
$\sum_Lm_L\le N$ and $\max_Lm_L\ll_\delta N^\delta$ for every $\delta>0$,
then the second-moment condition holds.
\end{lemma}
\begin{proof}
First consider nonnegative coefficients. Set $M=\sum_Lm_L$. Orthogonality gives
\[
\|Q_N\|_2^2=\sum_Lm_L^2\le(\max_Lm_L)\sum_Lm_L=(\max_Lm_L)M,
\]
This proves that the multiplicity bound suffices for the second-moment condition.
For a general $Q_N$ satisfying \eqref{eq:fair-L2}, for each fixed integer $k\ge2$ Bonami--Beckner hypercontractivity
\cite{Bonami,Beckner,BlaisTan} gives
\[
\|Q_N\|_{2k}\le(2k-1)^{d_N/2}\|Q_N\|_2.
\]
Since $k$ is fixed, the degree factor is $N^{o(1)}$. Choose the exponent in \eqref{eq:fair-L2} sufficiently small to obtain, for each fixed $k$ and $\eta>0$,
\[
\EE|Q_N|^{2k}\ll_{k,\eta}N^{k+\eta}.
\]
Markov's inequality then gives
\[
\PP\{|Q_N|>N^{1/2+\varepsilon}\}
\ll_{k,\eta}N^{-2k\varepsilon+\eta}.
\]
Choose $k$ sufficiently large and $\eta$ sufficiently small that $2k\varepsilon-\eta>2$. The bounds are then summable over all $N$.
The first Borel--Cantelli lemma gives, almost surely,
$Q_N=O_{\varepsilon,\mathcal X}(N^{1/2+\varepsilon})$.
Intersecting the events for $\varepsilon=1/j$ and for the finitely many statistics gives a single probability-$1$ event for every $\varepsilon>0$.
\end{proof}
\begin{remark}[Second moments and multiplicities]
The hypothesis of the lemma concerns the full second moment. The multiplicity bound is only a sufficient condition.
In the completely multiplicative one-point sum, the constant monomial comes from squares and has multiplicity about $\sqrt N$.
We verify the second-moment bound $O(N\log N)$ directly. The constant monomial is retained throughout. In the Borel--Cantelli argument,
summability of the tails allows dependence across different $N$.
\end{remark}

\subsubsection{Squarefree-supported model}
Let $g(n)=\mu^2(n)\prod_{q\mid n}\xi_q$ and use $U,V$ as defined in \eqref{eq:identity-sf}.
In $U$ distinct squarefree $a$ give distinct monomials, so $\|U\|_2^2\le N$.
For a nonzero term in $V$, write $a=du,b=dv$, where $d=(a,b)$ and $(u,v)=1$.
Squarefreeness gives $g(a)g(b)=\xi_{uv}$, and $d\mid N$.
For a fixed squarefree index $t=uv\le N^2$, the multiplicity of its monomial is at most
$\tau(N)\tau(t)\ll_\delta N^\delta$, where $\tau(n)$ counts divisors.
The total multiplicity is at most $N$. Lemma~\ref{lem:walsh} controls $U,V$ simultaneously, and \eqref{eq:identity-sf} gives the claim.
\begin{remark}[Coordinates in the squarefree-supported model]
The constant term with index $t=1$ comes from $a=b=N/2$ when that summand is squarefree. It is included in the multiplicities.
The argument uses only prime coordinates, so dependence among the higher prime-power coordinates has no effect in this model.
\end{remark}

\subsubsection{Fair model with independent prime powers}
The monomial index of an integer $n$ is
\[
L(n)=\{(q,v_q(n)):q\mid n\}.
\]
This index determines $n$ uniquely. Hence $\|S_f(N-1)\|_2^2=N-1$.
The product $f(a)f(b)$ has index $L=L(a)\mathbin\triangle L(b)$.
Fix $L$ and write $t=\prod_{(q,k)\in L}q^k\le ab\le N^2$.
Let $d$ be the product of the full prime powers whose exponents are equal and positive at the two summands, and write $a=du,b=dv$.
Then
\[
d\mid N,\qquad uv=t,\qquad d(u+v)=N.
\]
Indeed, all full prime powers with equal exponents at the two summands go into $d$. Prime powers with unequal positive exponents, or with a positive exponent on just one side, remain in $u,v$
and occur in the symmetric difference $L$, giving $uv=t$.
Thus for each fixed index $L$ the multiplicity satisfies
\[
m_L\le\tau(N)\tau(t)\ll_\delta N^\delta.
\]
The total multiplicity is $N-1$, so Lemma~\ref{lem:walsh} and \eqref{eq:identity} give the claim.
\begin{remark}[Cancellation of prime-power coordinates]
The integers $u,v$ can have a common factor. Different positive exponents of the same prime give different coordinates.
A coordinate cancels in the symmetric difference only if the exponents at both summands agree.
\end{remark}

\subsubsection{Common-sign prime-power model}
Now take $f(n)=\prod_{q\mid n}\xi_q$.
We first prove an elementary bound uniform in the prime support. For fixed $D,\delta>0$ and a squarefree integer $t\le N^D$,
\begin{equation}\label{eq:supported-count}
\Psi(N;t):=\#\{n\le N:\operatorname{rad}(n)\mid t\}
\ll_{D,\delta}N^\delta.
\end{equation}
Indeed, for any $s,\theta>0$, Rankin's bound gives
\[
\Psi(N;t)\le N^s\prod_{q\mid t}(1-q^{-s})^{-1}
\ll_{s,\theta}N^s t^\theta.
\]
The last step uses $(1-q^{-s})^{-1}\le q^\theta$ for sufficiently large $q$. The product over the finitely many small primes is bounded in terms of $s,\theta$ alone.
Taking $s=\delta/2$ and $\theta=\delta/(2D)$ gives \eqref{eq:supported-count}, with a constant independent of $t$ as it varies with $N$.

Group $S_f(N-1)$ by $t=\operatorname{rad}(n)$ to write it as
$\sum_t c_t\xi_t$. By \eqref{eq:supported-count},
$\max_t c_t\ll_\delta N^\delta$, while $\sum_t c_t=N-1$.
The one-point sum therefore satisfies the multiplicity condition in Lemma~\ref{lem:walsh}.

For $a+b=N$, the product index is
\[
t=\prod_{q:\,\mathbf1_{q\mid a}\ne\mathbf1_{q\mid b}}q\le ab\le N^2.
\]
Once $t$ is fixed, each prime dividing $a$ divides $Nt$. If it also divides $b$, it divides $N$. Otherwise it occurs in $t$.
Hence all corresponding $a$ satisfy
\[
\operatorname{rad}(a)\mid\operatorname{rad}(Nt),\qquad
\operatorname{rad}(Nt)\le N^3.
\]
By \eqref{eq:supported-count}, every monomial has multiplicity $O_\delta(N^\delta)$ uniformly, including the index $t=1$.
The total multiplicity is $N-1$, so the reflection sum also satisfies the multiplicity condition of Lemma~\ref{lem:walsh}. Together with \eqref{eq:identity} this proves the common-sign case.

\subsubsection{Fair completely multiplicative model}
Let $f(n)=\prod_q\xi_q^{v_q(n)}$.
Group the integers by the product $t$ of the primes with odd exponents. Each integer has a unique expression $n=tu^2$ with $t$ squarefree.
It follows that
\[
\|S_f(N-1)\|_2^2
\le\sum_{t\le N}\left\lfloor\sqrt{N/t}\right\rfloor^2
\ll N\log N.
\]
Squares contribute to the constant monomial. This bound includes the mean
$\EE S_f(N-1)=\lfloor\sqrt{N-1}\rfloor$.

To control the reflection sum we use the following representation bound.
\begin{lemma}[Uniform upper bound for quadratic representations]\label{lem:quadratic}
Let $A,B$ be coprime positive squarefree integers and $M$ a positive integer. Then
\begin{equation}\label{eq:quadratic-count}
\#\{x,y\ge1:Ax^2+By^2=M\}\le6\tau(AM).
\end{equation}
\end{lemma}
\begin{proof}
Take the imaginary quadratic field $K=\mathbb Q(\sqrt{-AB})$.
Each solution gives a distinct algebraic integer in its full ring of integers,
namely $\alpha=Ax+y\sqrt{-AB}$, with norm $AM$.
There are $e+1$ integral ideals of norm $q^e$ if $q$ splits, and at most $1$ if it is inert or ramified.
Unique factorization of ideals thus bounds the number of integral ideals of norm $AM$ by $\tau(AM)$.
Generators of a given principal ideal differ by units, and an imaginary quadratic field has at most $6$ units. This proves \eqref{eq:quadratic-count}.
The argument uses unique factorization of integral ideals in the Dedekind domain $\mathcal O_K$
\cite[Chapters~3--5]{Milne}, so the constant $6$ is independent of the discriminant and class number.
\end{proof}

The monomial index of the reflection product is the product $t$ of the primes with odd exponents in $ab$.
Write $a=du,b=dv$, with $d=(a,b)$ and $(u,v)=1$. Then
\[
u=t_1x^2,\quad v=t_2y^2,\quad
t_1t_2=t,\quad d\mid N,\quad
t_1x^2+t_2y^2=N/d.
\]
Here $t_1,t_2$ are coprime and squarefree.
By \eqref{eq:quadratic-count}, after fixing $d$ and $t_1\mid t$ there are at most $6\tau(t_1N/d)$ pairs.
Thus the multiplicity of each index satisfies
\begin{equation}\label{eq:cm-label-multiplicity}
m_t\le6\sum_{d\mid N}\sum_{t_1\mid t}\tau(t_1N/d)
\ll_\delta N^\delta,
\end{equation}
Since $t\le N^2$ and $t_1N/d\le N^3$, one can choose a sufficiently small fixed exponent in each of the three divisor bounds.
This includes $t=1$, and the total multiplicity is $N-1$. Apply Lemma~\ref{lem:walsh} to the reflection multiplicity bound
and to the second-moment bound $O(N\log N)$ for the one-point sum, then use \eqref{eq:identity} to obtain the bound with a positive power loss for the completely multiplicative model.

\subsubsection{Near-fair laws}
\begin{lemma}[Absolute continuity for square-summable biases]\label{lem:transfer}
Enumerate the independent basic signs as $\xi_j$. Let $\PP_0$ be the fair product measure and let $\PP_\delta$ be
the product measure with $\PP_\delta(\xi_j=-1)=1/2+\delta_j$, where $|\delta_j|\le1/2$ and $\sum_j\delta_j^2<\infty$.
Then $\PP_\delta\ll\PP_0$. Thus every probability-$1$ event for the fair measure also has probability $1$ for the biased measure.
\end{lemma}
\begin{proof}
On the fair probability space define the likelihood ratio for finitely many coordinates by
\[
L_J=\prod_{j=1}^J(1-2\delta_j\xi_j).
\]
This is a nonnegative martingale of mean $1$, and
\[
\EE_0 L_J^2=\prod_{j=1}^J(1+4\delta_j^2)
\le\exp\left(4\sum_j\delta_j^2\right)<\infty.
\]
Therefore $L_J$ converges in $L^2$ and $L^1$ to some $L\ge0$ with $\EE_0L=1$.
If $B$ depends only on the first $m$ coordinates, then for $J\ge m$ we have
$\EE_0(\ind_B L_J)=\PP_\delta(B)$. Taking the limit and using uniqueness of the product measure gives
\[
\PP_\delta(B)=\EE_0(\ind_B L)\quad\hbox{for every measurable }B,
\qquad \PP_\delta\ll\PP_0.
\]
This proves the required absolute continuity and the transfer of probability-$1$ events.
\end{proof}
In each model, intersect the events under the fair law for $\varepsilon=1/j$ and apply Lemma~\ref{lem:transfer}
to obtain a probability-$1$ event under the biased law. For arbitrary $\varepsilon>0$, choose a smaller $1/j$ to obtain the power bound.
This proves Theorem~\ref{thm:main}\textup{(ii)}.
\begin{remark}[Direction and scope of absolute continuity]
Lemma~\ref{lem:transfer} gives a sufficient condition for the infinite product-measure argument. For the classical background see Kakutani \cite{Kakutani}.
The required direction is $\PP_\delta\ll\PP_0$, which holds even when finitely many coordinates have probability $0$ or $1$,
whereas equivalence of the two measures can fail. For a fixed $\rho\ne1/2$ the sum of squared biases diverges, so Appendix~\ref{app:biased}
gives a separate argument. Divergence of the sum alone does not determine the size of the error.
\end{remark}

\section{Proof of the exact extremes in the completely multiplicative model}\label{sec:cm-proof}
\paragraph{Proof outline}
We first bound the limit superior of the normalized absolute reflection sum by $2$.
We then show that each specified short interval contains positive and negative values approaching this bound.
Finally we use a known almost-sure bound for the one-point sum to pass from the reflection sum to the four representation counts.

For the upper bound, we discard a small set of positions and split the remaining labels according to whether their largest prime factor exceeds $Y_0$.
A predictable square sum controls the part with large prime factors. We group the rest by decreasing prime prefixes,
estimate the groups at which the prefix procedure stops, and work back through the levels to the original sum.
For the lower bound, we estimate the probability of a large value at one candidate target and then the probability of simultaneous large values at two targets.
These estimates show that at least one target in the short interval reaches the threshold.
Write $L=\log N$, $l=\log L$, $U=(Ll)^{1/3}$, and $D=L/U$.
Table~\ref{tab:cm-constants} lists the cutoffs used for the upper bound and the corresponding probability requirements.

\begin{table}[htbp]
\centering\small
\begin{tabular}{@{}>{\raggedright\arraybackslash}p{.38\linewidth}>{\raggedright\arraybackslash}p{.57\linewidth}@{}}
\hline
Cutoff or parameter & Role in the proof\\
\hline
$Z=N^{1/10}$, with summands at distance at least $N^{9/10}$ from the boundary &
The retained kernels are at least $N^{17/10}/2$. The failure probability for the discarded sum at a fixed target is $O(N^{-2})$.\\
$Y_0=e^{9D}$,  and $M_0=Ne^{-Ul/40}$ &
In the large-prime range take $k\sim U/2$. The identical-pair term gives an exponent margin of $9-12(1/2)^2=6$.
The mass saving in the small-prime range covers the $\exp(O(U))$ loss incurred in recovering the original sum.\\
Stop when the prefix product reaches $e^{6.9D}$. Take $k\sim0.44U$. &
The tail exponent at a fixed target is $(0.44/2)(6.9-12(0.44)^2)=1.006896>1.005$.\\
Stop when the largest prime factor of the remainder is at most $e^{0.30U}$. Take $k\sim U/6$. &
The tail exponent at a fixed target is $(20-e^3/3)/12>1.005$.
The prefix is at most $e^{15.9D}<N^{1/4}$, and the depth is at most $\lceil23D/U\rceil+1$.\\
\hline
\end{tabular}
\caption{Cutoffs and summability conditions for the completely multiplicative reflection sum.}\label{tab:cm-constants}
\end{table}

\subsection{An upper bound for the reflection sum at all integer targets}\label{sec:cm-upper}
We first prove Theorem~\ref{thm:cm-extrema} for fair basic signs, keeping its notation $s_N$. Fix a sufficiently large $N$ and note that $U^2=Dl$.
\paragraph{Discarding positions and grouping Walsh labels}
Write $Q_N(a)=a(N-a)=t_a z_a^2$ with $t_a$ squarefree, and set $\xi_t=\prod_{p\mid t}\xi_p$.
All positions are taken over $1\le a<N$, including the midpoint when the target is even.

\begin{lemma}\label{lem:cm-thin}
There is a deterministic set of positions $E_N$ that can be discarded so that
$\PP(|\sum_{a\in E_N}\xi_{t_a}|>\sqrt N)\ll N^{-2}$.
Among the remaining labels, the multiplicities $M_t$ are at most $2$. For ordered quadruples of pairwise distinct labels whose product is a square,
the total weighted count, with weight equal to the product of their four multiplicities, satisfies $B_4\ll N^{5/3}$.
\end{lemma}
\begin{proof}
Set $Z=N^{1/10}$. First discard positions for which $\min(a,N-a)<N^{9/10}$ or $z_a>Z$.
If $d=(a,N),\ a=du,\ N-a=dv$, then $(u,v)=1$ and $z_a=d\,s(u)s(v)$,
where $s(n)^2$ is the largest square divisor of $n$. Using
$\sum_{n\le y}s(n)\le\sum_{k\le\sqrt y}k\lfloor y/k^2\rfloor\ll y\log(2y)$,
Cauchy--Schwarz for the sum over $u+v=N/d$ gives
$\sum_a\sqrt{z_a}\ll N\log(2N)\sum_{d\mid N}d^{-1/2}$.
Thus at most $N^{19/20+o(1)}$ positions have been discarded.

Next discard positions with a remaining equal-kernel partner other than their reflection, that is, $t_a=t_b$ with $b\ne a,N-a$.
Then $z_a\ne z_b$ and the ratio satisfies $Q_N(a)/Q_N(b)=z_a^2/z_b^2\ne1$.
We use Pylaiev's fixed-ratio counting lemma \cite[Lemma~4.4]{Pylaiev} twice below.
Fix a positive rational ratio $\varrho\ne1$ whose reduced numerator and denominator are at most $N^6$.
The number of ordered solutions to $Q_N(a)/Q_N(b)=\varrho$ is $N^{o(1)}$, uniformly in these ratios.
In applying the lemma, reflect each variable separately into the half-interval up to $N/2$, at a cost of at most a factor of four.
If one variable equals $N/2$, the other has at most two values for a fixed ratio.
Since $z_a,z_b\le Z$, this discards at most another $Z^2N^{o(1)}$ positions, leaving $M_t\le2$.
The maximum multiplicity before any discarding is uniformly $N^{o(1)}$ by \eqref{eq:cm-label-multiplicity}.
Thus the second moment of the discarded sum is at most $|E_N|N^{o(1)}$, and hypercontractivity at the fixed moment order $100$ costs a factor
$N^{o(1)}$ in the norm, since the degree is $O(L/l)$. Markov's inequality at $\sqrt N$ gives
a tail probability of $N^{-5/2+o(1)}$.

For a remaining quadruple, fix the first position $a$. Its label satisfies $t_a\ge N^{17/10}/2$.
Each of its prime factors occurs in at least one of the other three labels, so some $t_b$ satisfies
$G=(t_a,t_b)\ge t_a^{1/3}$. Enumerate $G\mid t_a$. Since $G\mid b(N-b)$,
there are at most two roots, $0,N$, modulo each $p\mid G$. For $b$, the Chinese remainder theorem gives at most
$N^{13/30+o(1)}$ choices. Let $T$ be the squarefree kernel of $t_at_b$. Then $T>1$.
For the other two positions $c,e$, remove the greatest common divisor of $t_c,t_e$. The two coprime factors that remain have product $T$.
Enumerating their possible factorizations and the values $z_c,z_e\le Z$ fixes the ratio $Q_N(c)/Q_N(e)$.
Its square class is $T$ and its numerator and denominator are at most $N^6$. Applying the ratio bound again gives
$B_4\le N^{1+13/30+1/5+o(1)}\ll N^{5/3}$.
\end{proof}

The discarding lemma bounds the remaining multiplicities and the count of four-label relations. Both estimates are used to control fluctuations of the group square sum.
From now on count only the remaining labels. Partition them into disjoint groups, each with a specified squarefree prefix $g\mid t$.
Set $A_g=\sum_{t\text{ in }g} M_t\xi_{t/g}$, $W=\sum_g A_g^2$, and
$b_g(v)=\sum_{t\text{ in }g}M_t^2v^{\omega(t/g)}$.
Then $\EE W=\sum_g b_g(1)\le2N$. For $k\sim \alpha_0 U$ with $\alpha_0>0$ and for $v=k-1$,
\begin{equation}\label{eq:cm-group}
 \|W-\EE W\|_k^2\le 2\sum_g b_g(v)^2+N^{2-1/12+o(1)} .
\end{equation}
For a Walsh polynomial $\sum c_s\xi_s$, the noise operator multiplies each coefficient by
$(k-1)^{-\omega(s)/2}$. Thus $2\to k$ hypercontractivity bounds the squared norm by
$\sum_s c_s^2 v^{\omega(s)}$. In $W-\EE W$, each nonempty label $s$ is the symmetric difference of two
distinct labels in one group. If two ordered pairs with the same difference share a label, they are the same pair or reversals.
Their weighted contribution is at most $2\sum_g b_g(v)^2$. Otherwise the four full labels $t_i$ are distinct and their product is a square.
We also have $2\omega(s)\le\sum_i\omega(t_i)$. Apply H\"older's inequality to this part with $h=4/3$
to bound its contribution by
\[
 B_4^{1/4}\left[\EE\left(\sum_t M_t v^{h\omega(t)/2}\xi_t\right)^4\right]^{3/4}.
\]
By $2\to4$ hypercontractivity, the expression in square brackets is at most $(\sum_t M_t^2(3v^h)^{\omega(t)})^2$.
For $w=3v^h$, use $M_t\le2$, reflection symmetry, and Cauchy--Schwarz to obtain
\[
 \sum_t M_t^2 w^{\omega(t)}
 \le 2\sum_{n<N}w^{2\omega(n)}
 \le 2N\prod_{p\le N}(1+w^2/p)=N^{1+o(1)}
\]
Here $w^2 l=o(L)$, which proves \eqref{eq:cm-group}. The count includes terms from different groups.

\begin{lemma}[Local weight of a prefix]\label{lem:cm-prefix}
Suppose that $g\le N^{1/4}$, $\omega(g)=O(D/U)$, and $v\sim \alpha_0 U$.
Uniformly over the groups above,
\[
 b_g(v)\le N g^{-1+o(1)}\exp((6\alpha_0^2+o(1))D).
\]
If also $g\le e^{6.9D}$ and every prime factor of each remainder $t/g$ is at most
$e^{0.30U}$, then
\[
 b_g(v)\le (N/g)\exp((6\alpha_0^2 e^3-10+o(1))D).
\]
\end{lemma}
See Appendix~\ref{app:cm-technical} for the proof. We use this lemma to bound the predictable square sum for the large primes and the prefix-group weights at the stopping points for the small primes.

\noindent\textbf{Large primes.}
Reveal the signs in increasing prime order. For any deterministic subset of labels, group by the largest prime factor $p$. The coefficient $A_p$
depends only on signs at primes below $p$. Put $Y_0=e^{9D}$. For $Y_0<p\le N^{1/4}$, take $k\sim U/2$.
Lemma~\ref{lem:cm-prefix} and $\sum_{p>Y_0}p^{-2+o(1)}\le e^{(-9+o(1))D}$ give
$\sum_p b_p(v)^2\le N^2e^{(-6+o(1))D}$. By \eqref{eq:cm-group},
for any fixed $\eta>0$ the probability that $W$ for this range differs from its mean by more than $\eta N/2$ is at most $N^{-3/2+o_\eta(1)}$,
because $kD\sim L/2$. The same applies to any deterministic subrange, since removing groups can only decrease the collision sum for nonempty Walsh coefficients.
In the range $p>N^{1/4}$, we have $b_p(1)\le4N/p$ and $\sum_p b_p(1)\le2N$.
Separating repeated pairs from quadruples of distinct labels gives $\operatorname{Var}(W)\le2\sum_p b_p(1)^2+B_4\ll N^{7/4}$.
The centered polynomial has degree $O(L/l)$. Hypercontractivity at the fixed order $32$ gives a probability of $N^{-4+o(1)}$ at the same threshold.
The sum of the two means is at most $2N$. Thus for $V_N=\sum_{p>Y_0}A_p^2$ we have
\[
 \PP(V_N>(2+\eta)N)\le N^{-3/2+o_\eta(1)}+N^{-4+o_\eta(1)} .
\]

We use the following martingale bound for predictable Euclidean vectors $V_p$.
\begin{equation}\label{eq:cm-vector}
 \PP\left(\left\|\sum_p\xi_p V_p\right\|\ge t,\ \sum_p\|V_p\|^2\le b^2\right)
 \le 2e^{-t^2/(2b^2)} .
\end{equation}
Indeed, average the power series of $\cosh(\lambda\|x\pm y\|)$ over the two signs. Each term is nondecreasing in
$|\langle x,y\rangle|$, so the average is bounded by its value when the vectors point in the same direction, namely
$\cosh(\lambda\|x\|)\cosh(\lambda\|y\|)$. Iterating conditional expectations with $\cosh s\le e^{s^2/2}$
gives
\[
 \EE\!\left[\cosh\!\left(\lambda\left\|\sum_p\xi_p V_p\right\|\right)
       \exp\!\left(-\tfrac12\lambda^2\sum_p\|V_p\|^2\right)\right]\le1.
\]
Taking $\lambda=t/b^2$ proves \eqref{eq:cm-vector}.
The joint probability that $V_N\le(2+\eta)N$ and that the absolute value of the large-prime contribution reaches $(2+\eta)s_N$
is at most $2N^{-(2+\eta)^2/[2(2+\eta)]}=2N^{-1-\eta/2}$.

\noindent\textbf{Small primes and stopping groups.}
For the remaining labels put $M_0=N\exp(-Ul/40)$. We first prove
\[
 \sum_{t:\,p\mid t\,\Rightarrow\,p\le Y_0}M_t^2\le M_0 .
\]
For $t\le N$, count smooth integers using $M_t^2\le4$.
For $t>N$, let $d\le N$ be the longest increasing prime prefix of $t$ that does not exceed the target. Then $N/Y_0<d\le N$.
Since $d\mid Q_N(a)$, the number of positions for a given $d$ is $\le2(N/d)2^{\omega(d)}$, and $M_t^2\le2M_t$.
Put $\lambda=l/(4\log Y_0)$. Rankin's bounds for the two parts are, respectively,
\[
 4N^{1-\lambda}\prod_{p\le Y_0}(1+p^{\lambda-1}),\qquad
 4N(Y_0/N)^\lambda\prod_{p\le Y_0}(1+2p^{\lambda-1}).
\]
Convexity gives $p^\lambda-1\le(e^{l/4}-1)\log p/\log Y_0$, so the logarithms of the two products are
$O(l+e^{l/4})=o(Ul)$. The desired mass bound follows from $\lambda L=Ul/36$.

Starting from $g=1$, stop if $g\ge e^{6.9D}$ and call this a large-product terminal group.
Otherwise, stop if the remainder is $1$ or if its largest prime factor is $\le e^{0.30U}$, and call this a smooth terminal group.
Send all other labels to $gp$, where $p$ is the largest prime factor of the remainder. Some labels can stop at a node while others continue.
The depth is at most $h=\lceil23D/U\rceil+1$, and all prefixes satisfy $g\le e^{15.9D}<N^{1/4}$.
For each of the two terminal families its square sum $W$ satisfies $\EE W\le M_0$.
Take $k\sim0.44U$ for the large-product terminals and $k\sim U/6$ for the smooth terminals.
Square the bound in Lemma~\ref{lem:cm-prefix} and sum, using $\sum_{g\ge e^{6.9D}}g^{-2+o(1)}
\le e^{(-6.9+o(1))D}$ for the former and $\sum_g g^{-2}=O(1)$ for the latter. Substitute these estimates into
\eqref{eq:cm-group} and apply Markov's inequality to the centered square sum at $M_0$.
As $Ul=o(D)$ and $UD=L$, the probabilities that the two square sums exceed $2M_0$ are bounded, respectively, by
\[
 N^{-(0.44/2)(6.9-12(0.44)^2)+o(1)},\qquad
 N^{-(20-e^3/3)/12+o(1)} .
\]
The first exponent constant is $(11/50)(69/10-12(11/25)^2)=1.006896$.
By $e^3<21$, the second constant is greater than $13/12$.
For either fixed ratio of moment order to scale, $o(1)$ tends to zero as $N\to\infty$, uniformly over deterministic subfamilies of terminals.
Thus both failure probabilities are eventually $O(N^{-1.005})$.

Index the vector at level $j$ by the prefixes $g$ reached at that level. Let $Z_j$ be the norm of the vector of sums over all their arriving remainders.
If both terminal families satisfy their bounds, the sum of the norms of the two parts stopping at the current level is $\le3\sqrt{M_0}$.
The part that continues is $\sum_p\xi_p V_p$, where $(V_p)_g$ is the remainder sum at the child node $gp$
and only uses primes below $p$. A decreasing prefix has a unique parent, so $\sum_p\|V_p\|^2=Z_{j+1}^2$.
Define $b_h=3\sqrt{M_0}$ and $b_j=3\sqrt{M_0}+\sqrt{6L}\,b_{j+1}$.
If the terminal bounds hold but $Z_0>b_0$, then some level exceeds $b_j$ while the next level stays within its bound.
By \eqref{eq:cm-vector} and a union bound, this has probability at most $2hN^{-3}$. Also $h\log(6L)=O(U)$, so
$b_0\le6\sqrt{M_0}(\sqrt{6L})^h=o(\sqrt N)$.
Thus outside an event of probability $O(N^{-1.004})$, the sum of retained labels with $P^+(t)\le Y_0$,
including any constant term, has absolute value at most $b_0$. Restore the sum discarded in Lemma~\ref{lem:cm-thin}.
The first Borel--Cantelli lemma and a countable intersection over $\eta=1/m$ give $\limsup_N |C_f(N)|/s_N\le2$.

\subsection{Two-sided block extremes of the reflection sum}\label{sec:cm-lower}
Fix $0<h<1/10$ and put $\delta=h/4$ and $\nu=h/64$. Along $X=2^j$ take $y=X^\delta$ and the candidate set
$T=[X,X+\lfloor X^{1-h}\rfloor)\cap\mathbb Z$. For each $i\in[X,2X)$,
apply the discarding rule of Lemma~\ref{lem:cm-thin} and write $E_i$ for the discarded sum. The rule is symmetric under $a\leftrightarrow i-a$,
and the midpoint of an even target, with $z=i/2$, is discarded. Thus each retained label occurs at exactly two positions and satisfies
$z\le i^{1/10}$,  and $t\ge i^{17/10}/2$.

\begin{lemma}\label{lem:cm-cross}
Consider any four specified targets in $[X,2X)$, keeping their positions according to the rule above.
The number of ordered position quadruples with pairwise distinct labels whose product is a square is $\ll X^{5/3}$, uniformly over the targets.
For distinct $i,k\in[X,2X)$, the number of pairs of retained positions with the same label is at most $2\gcd(i,k)+X^{1/5+o(1)}$.
\end{lemma}
\begin{proof}
Fix the first position. Each prime factor of its label $t$ occurs at least once more,
so one of the other labels has greatest common divisor $G\ge t^{1/3}\gg X^{17/30}$ with it.
Enumerate the partner position among the three slots and the divisors $G\mid t$. For its target $k$, use the roots $0,k$ modulo each $p\mid G$
to obtain at most $X^{13/30+o(1)}$ possible partner positions. Fixing this pair fixes the squarefree kernel $D_0>1$ of the product of their labels.
Removing the greatest common divisor from the labels of the other two positions $c,e$ leaves $A,B$ with $AB=D_0$.
Enumerate the $\tau(D_0)$ possible factorizations and the square parts $z,w\le(2X)^{1/10}$ for the two positions.
This fixes the reduced ratio $s/t'=Az^2/(Bw^2)$ of their two summand products. The positive rational ratio has square class $D_0$.
If the two targets are $i,k$, put $x=2c-i,\ u=2e-k$. Then
\[
 t'x^2-su^2=t'i^2-sk^2\ne0 .
\]
A zero right-hand side would force $s/t'=(i/k)^2$. In the real quadratic field $\mathbb Q(\sqrt{st'})$,
the pairs of positions inject into algebraic integers $t'x+u\sqrt{st'}$ with fixed norm
$K=t'(t'i^2-sk^2)$. Since $s,t',|K|\le X^{O(1)}$, the absolute values under both real embeddings lie in
$[X^{-O(1)},X^{O(1)}]$. There are at most $\tau(|K|)$ integral ideals of norm $|K|$,
by counting prime-ideal factors according to the splitting type of each rational prime \cite[Ch.~3]{Milne}.
If $\alpha,\beta$ are two eligible generators of the same ideal, their ratio $\alpha/\beta$ is a unit of
norm $K/K=1$, since the equation fixes the signed norm $K$.
Fix one real embedding and take the ratio $\lambda>1$ of the larger to the smaller of two distinct absolute values.
This is the value of the norm-$1$ unit $\pm\alpha/\beta$ (or its inverse) at that embedding. Its value at the other embedding is $\lambda^{-1}$.
Its integral trace satisfies $\lambda+\lambda^{-1}\ge3$, so $\lambda\ge(3+\sqrt5)/2$ uniformly over the fields.
In the polynomial-size interval above, each ideal therefore has $O(\log X)$
eligible generators, counting both signs. There are $X^{o(1)}$ completions for each ratio, so the total is at most
$X^{1+13/30+1/5+o(1)}\ll X^{5/3}$.

For the second assertion, equality of labels gives $w^2 a(i-a)=z^2 b(k-b)$.
If $wi=zk$, then $a/i=b/k$ or $a/i=1-b/k$, and the two lines contain at most
$2\gcd(i,k)\le2|i-k|$ interior lattice points in total. Otherwise
\[
 [w(2a-i)-z(2b-k)][w(2a-i)+z(2b-k)]=w^2 i^2-z^2 k^2\ne0 .
\]
For each $z,w$ there are $X^{o(1)}$ signed factor pairs, each determining $a,b$. This gives $X^{1/5+o(1)}$ in total.
\end{proof}

\noindent\textbf{The frozen and subsequently revealed parts.}
First reveal the signs of primes $p\le y$. Write $C_f(i)=E_i+L_i+H_i$, where $L_i$ is the sum over retained labels
containing only primes $p\le y$. The remaining part is
\[
 H_i=\sum_{p>y}\xi_p A_{p,i},\qquad
 A_{p,i}=\sum_{P^+(t)=p}M_t(i)\xi_{t/p},\qquad
 b_{p,i}=\sum_{P^+(t)=p}M_t(i)^2 .
\]
Continue revealing signs in increasing prime order, so $A_{p,i}$ is known before $\xi_p$ is revealed. By orthogonality,
$\sigma_i^2:=\EE H_i^2=\sum_{p>y} b_{p,i}$, and uniformly
\begin{equation}\label{eq:cm-revealed}
 2(1-\delta-o_h(1))i\le\sigma_i^2\le 2i .
\end{equation}
Indeed, after discarding $a<i/\log i$ the remaining positions satisfy $\log a\ge\log i-\log\log i$.
The sum of $\log a$ over all $y$-smooth positions is at most
$i\sum_{p\le y}\sum_{k\ge1}(\log p)/p^k=i(\log y+O(1))$, so there are at most $(\delta+o_h(1))i$ such positions.
Discard also positions for which some $p>y$ satisfies $p^2\mid a$ or $p\mid(a,i)$. There are $O_h(i/y)$ such positions,
because $\sum_{m>y}m^{-2}\ll1/y$ and $i$ has at most $O_h(1)$ such large prime factors.
Each remaining position has a prime $p>y$ that occurs exactly once in $a(i-a)$. The discarding rule removes only $o(i)$ positions,
and each retained label satisfies $M_t(i)^2=2M_t(i)$. This gives \eqref{eq:cm-revealed}.

Let $Y_i$ be the cutoff $Y_0$ from Section~\ref{sec:cm-upper} for $N=i$.
For a fixed $h$, eventually $Y_i<y<i^{1/4}$ throughout the block.
The predictable square sum $W_i$ for the frozen part on $(Y_i,y]$ satisfies
$\EE W_i\le2i-\sigma_i^2\le(2\delta+o_h(1))i$.
The preceding deterministic-subrange estimate gives $\PP(W_i>3\delta i)\le i^{-3/2+o_h(1)}$.
The joint probability that $W_i\le3\delta i$ and that the absolute value of the subrange sum exceeds $\sqrt{8\delta}\,s_i$
is at most $2i^{-4/3}$ by \eqref{eq:cm-vector}. Include the exceptional probabilities for the $o(\sqrt i)$ bound below $Y_i$
and for $|E_i|\le\sqrt i$. A union bound over at most $X$ targets gives
\begin{equation}\label{eq:cm-freeze}
 \PP\bigl(\exists i\in T:\ |E_i+L_i|>3\sqrt\delta\,s_i\bigr)\ll_h X^{-1/1000}.
\end{equation}

\noindent\textbf{Predictable covariance.}
Prescribe a direction $e\in\{-1,1\}$ and take $d=1$ or $2$ distinct candidate targets, indexed by $i$.
Set $a_p=(eA_{p,i}/\sigma_i)_i$ and $\Gamma=\sum_{p>y}a_pa_p^{\mathsf T}$.
Uniformly over the selected coordinates we have
\begin{equation}\label{eq:cm-good}
 \PP(\mathcal G^c)\ll_h X^{-80},\qquad
 \mathcal G=\{\max_{p>y}|a_p|\le X^{-\nu},\
                  \|\Gamma-I_d\|_{\rm op}\le X^{-\nu}\}.
\end{equation}
To prove this, write $Q_{ik}=\sum_{p>y}A_{p,i}A_{p,k}$.
Orthogonality gives $\EE Q_{ii}=\sigma_i^2$. For distinct targets Lemma~\ref{lem:cm-cross} bounds the mean
by $2\gcd(i,k)+X^{1/5+o(1)}\ll X^{1-h}$, since $\gcd(i,k)\le|i-k|$. Also,
\[
 \operatorname{Var}(Q_{ik})\le2\sum_{p>y}b_{p,i}b_{p,k}+O(X^{5/3}).
\]
The two labels in each row share a largest prime, whose cancellation preserves the symmetric difference. Centering removes the empty difference.
If two rows have the same nonempty difference and share a label, the unordered pairs coincide. Pairs in the same order contribute at most
$\sum_p b_{p,i}b_{p,k}$, and Cauchy--Schwarz gives the same bound for reversed pairs. For four distinct labels use Lemma~\ref{lem:cm-cross}.
Since $p\mid a(i-a)$, we have $b_{p,i}\le4i/p$. Combining this with $\sum_{p>y}b_{p,i}\le2i$,
the $L^2$ norms of the normalized, centered matrix entries and of each $A_{p,i}/\sigma_i$ are $\ll_h X^{-\delta/2}$.
Their Walsh degrees are $O(\log X/\log\log X)$. Take a fixed even moment order $\kappa>84/\nu$.
Hypercontractivity costs only a factor $X^{o_\kappa(1)}$ in the norm. Markov's inequality and a union bound over at most $2X$ primes and
at most four matrix entries give \eqref{eq:cm-good}.

\begin{lemma}[Relative moderate deviations]\label{lem:cm-relative}
Put $z=\sqrt{2(1-2h)\log X}$ and $u_X=\PP(G\ge z)$, where $G$ is standard normal.
Under the fair law, uniformly over the selection of one or two targets,
\[
 \PP(eH_i\ge z\sigma_i\ \text{for all selected }i)
   =u_X^d\bigl(1+O_h(X^{-\nu/100})\bigr).
\]
\end{lemma}
See Appendix~\ref{app:cm-technical} for the proof. The one-target estimate controls the mean number of exceedances, and the two-target estimate controls its second moment.

\noindent\textbf{Positive and negative extremes in each block.}
For a fixed direction write $Z_e=\sum_{i\in T}\ind_{\{eH_i\ge z\sigma_i\}}$ and $n=|T|$.
Lemma~\ref{lem:cm-relative} gives $\EE Z_e=nu_X(1+O_h(X^{-\nu/100}))$,
$\EE[Z_e(Z_e-1)]=n(n-1)u_X^2(1+O_h(X^{-\nu/100}))$.
Since $nu_X\asymp_h X^h/\sqrt{\log X}$, Chebyshev's inequality gives
$\PP(Z_e=0)\ll_h X^{-h/2}+X^{-\nu/100}$. Take a union bound over the two directions and the exceptional event in \eqref{eq:cm-freeze}.
The resulting probabilities are summable along $X=2^j$. In every other block, for each direction some $i\in T$ satisfies
\[
 \frac{eC_f(i)}{s_i}\ge
 2\sqrt{(1-\delta-o_h(1))(1-2h)\frac{\log X}{\log i}}-3\sqrt\delta .
\]
Apply the first Borel--Cantelli lemma for each fixed $h$, then intersect over $h=1/m$ with $m\ge11$.
Since $T\subset I_X$ eventually and $\log i/\log X\to1$ uniformly, let $m\to\infty$.
The lower limit of the block maxima for either direction is at least $2$. The all-integer upper bound then gives both block limits and the limit superior and limit inferior along the full sequence.

\noindent\textbf{The one-point sum and the four counts.}\par
In the fair completely multiplicative model, $f$ is an extended Rademacher
random multiplicative function in the sense of Atherfold. Its independent fair
prime signs are extended completely multiplicatively to all positive integers.
Taking $\varepsilon=1/4$ in Atherfold's theorem \cite[Theorem~1]{AtherfoldExtended}
gives, almost surely for all sufficiently large $x$,
\[
 |S_f(x)|\ll_{\mathcal X}\sqrt{x}(\log\log x)^{1/2}
 =o(\sqrt{x\log x}).
\]
Intersect this event with the probability-one event for the reflection sum established above.
By \eqref{eq:identity},
\[
 \begin{aligned}
 D_{\sigma\tau}(N)
 &=\frac{\sigma\tau C_f(N)}{4s_N}
   +\frac{(\sigma+\tau)S_f(N-1)-1}{4s_N}\\
 &=\frac{\sigma\tau C_f(N)}{4s_N}
   +O_{\mathcal X}\!\left(\frac{(\log\log N)^{1/2}}{\sqrt{\log N}}\right).
 \end{aligned}
\]
The remainder tends to zero uniformly on $I_X$. For opposite signs the positive and negative extremes are interchanged. The two constants for each of the four counts are $1/2,-1/2$. Taking the maximum of the four absolute errors gives $\limsup_N\Delta_f(N)/s_N=1/2$.

In the completely multiplicative model the basic coordinates are the prime signs. Intersect the events for the all-integer upper bound and the two-sided short-interval extremes,
then use the one-way absolute continuity $\PP_\delta\ll\PP_0$ of Lemma~\ref{lem:transfer}.
This gives the same pathwise conclusions when the prime signs are independent and their biases are square summable, by transferring the common probability-$1$ event for a fixed law.
The exact constants come from the all-integer upper bound, the block lower bound, and the one-point estimate in this section.
This completes the proof of Theorem~\ref{thm:cm-extrema}.

\section{Density, scope, and open problems}\label{sec:scope}
\subsection{Zero-density negative primes and an explicit threshold}\label{sec:sparse}
We use a concrete family of biases to show how a threshold arises. For every $\alpha>0$, the negative primes have relative density zero.
When $\alpha\le1$, their probability weight is still large enough to make the error $o(N)$.
When $\alpha>1$, truncation to finitely many primes and the Chinese remainder theorem show that the mean sign on the integers is positive.
The marginal count then rules out the quarter law.

\medskip\noindent\textbf{Step 1. Density of negative primes.}
Let $q_j$ be the prime of index $j$ and put $I_j=\ind_{\{X_{q_j,1}=-1\}}$.
Since $\rho_{q_j}\to0$, the average of the expectations tends to $0$. Hoeffding's inequality \cite{Hoeffding} gives
\[
\PP\left\{\left|J^{-1}\sum_{j\le J}(I_j-\rho_{q_j})\right|>\eps\right\}
\le2e^{-2\eps^2J}.
\]
Summability and Borel--Cantelli imply that almost surely
\[
\frac{\#\{q\le x:f(q)=-1\}}{\pi(x)}\to0.
\]
In Step 3 the density of signs on the integers is computed separately using a multiplicative mean.

\medskip\noindent\textbf{Step 2. The scale in the divergent case.}
Put $T=\log\log N$. Mertens' second theorem \eqref{eq:mertens-reciprocal} and partial summation give
\begin{equation}\label{eq:integral}
 B(N)=\int_1^T t^{-\alpha}\,dt+O_\alpha(1)\qquad(0<\alpha\le1).
\end{equation}
Indeed, the bounded remainder in $\sum_{q\le x}1/q=\log\log x+O(1)$ and the total variation of the monotone weight $t^{-\alpha}$ produce only an
$O_\alpha(1)$ error. The shift to $q+e^e$ in \eqref{eq:family} also has a summable effect on the reciprocal-prime sum.
Since $\log\log y_N=\log(3T)$, for $0<\alpha<1$
equation \eqref{eq:integral} gives $B(y_N)\ll_\alpha 1+(\log(3T))^{1-\alpha}$,
whereas $B(N)\sim T^{1-\alpha}/(1-\alpha)$.
For $\alpha=1$, the same formula gives $B(y_N)=\log\log(3T)+O(1)$ and $B(N)=\log T+O(1)$.
In both cases $B(y_N)=o(B(N))$, so eventually $2A_N=2(B(N)-B(y_N))\ge B(N)$.
Substituting \eqref{eq:integral} into \eqref{eq:weight-all} gives $E_\alpha(N)$, with squarefree main term $F(N)/4$.

\medskip\noindent\textbf{Step 3. The positive-mean obstruction in the convergent case.}
If $\alpha>1$, then $\sum_q\rho_q/q<\infty$. The nonnegative random sum has finite expectation, so almost surely the set $\mathcal S$ of negative primes satisfies
\[
\sum_{q\in\mathcal S}1/q<\infty.
\]
Fix such a realization. Since every higher prime-power sign is $+1$, a negative prime $q$ flips the sign of $f(n)$ only when $v_q(n)=1$.
For a finite subset $E\subset\mathcal S$, the Chinese remainder theorem gives the mean over the integers in the finite model,
\[
M_E=\prod_{q\in E}(1-2/q+2/q^2).
\]
The limit superior of the absolute difference between the means in the infinite and finite models is at most $2\sum_{q\in\mathcal S\setminus E}1/q$, since only integers divisible by these omitted primes can change signs.
Hence
\begin{equation}\label{eq:positive-mean}
M(f):=\lim_{x\to\infty}\frac1x\sum_{n\le x}f(n)
=\prod_{q\in\mathcal S}(1-2/q+2/q^2)>0.
\end{equation}
Every factor is positive, and their deviations from $1$ are summable, so the product is strictly positive.
The density of negative integers is $(1-M(f))/2<1/2$.
On the other hand, the marginal counting identity gives
\[
\sum_{\tau\in\{-1,1\}}R^f_{-\tau}(N)=\#\{1\le a<N:f(a)=-1\}.
\]
If all four counts were $N/4+o(N)$, the right-hand side would be $N/2+o(N)$. This contradicts the density strictly below $1/2$.

\subsection{Marginal density and additive reflection}\label{sec:patterns}
\subsubsection{Pair distribution and the reflection obstruction}
The structural example at density one half below shows why marginal balance cannot replace a bound for the reflection correlation.
Fix a nonvanishing realization $f$ and use the local notation
\[
s_N=\frac{S_f(N-1)}{N-1},\qquad c_N=\frac{C_f(N)}{N-1}.
\]
By \eqref{eq:identity},
\[
\frac{R^f_{\sigma\tau}(N)}{N-1}
=\frac{1+(\sigma+\tau)s_N+\sigma\tau c_N}{4}.
\]
Thus all four counts satisfy $R^f_{\sigma\tau}(N)=N/4+o(N)$ simultaneously
if and only if $s_N,c_N\to0$. Marginal balance controls only $s_N$.
The reflection correlation $c_N$ requires a separate estimate.
More precisely,
\[
\max_{\sigma,\tau}\left|R^f_{\sigma\tau}(N)-\frac{N-1}{4}\right|
=\frac{2|S_f(N-1)|+|C_f(N)|}{4}.
\]
In the squarefree model, \eqref{eq:identity-sf} gives the corresponding relation, with center $F(N)/4$ and terms $U(N)$ and $V(N)$.
These identities separate the effects of the marginal bias and the reflection correlation. A bound for the former cannot replace a bound for the latter.
The next proposition shows that even if the relative density of negative primes and the natural density of negative integers both equal $1/2$,
structure modulo $3$ can force $c_N=-1$ along a target sequence, ruling out same-sign pairs there.

\begin{proposition}[A structural obstruction at density one half]\label{prop:mod3}
There exists a genuinely non-completely multiplicative sign function for which the relative density of negative primes and the natural density of negative integers are both $1/2$, but for all $k\ge1$ one has $R_{--}(3^k)=R_{++}(3^k)=0$.
\end{proposition}
\begin{proof}
Let $\chi_3$ be the nonprincipal real character modulo $3$. Write $n=3^t u$, with $3\nmid u$, and define
\[
c_0=1,\quad c_2=-1,\quad c_t=1\ (t\ne2),\qquad f(n)=c_t\chi_3(u).
\]
The function is multiplicative on coprime arguments. Since $f(3)=1,f(9)=-1$, it is not completely multiplicative.
At each fixed $3$-adic valuation, the two unit classes have equal density within that layer. Discarding $t\ge T$ and then letting $T\to\infty$ shows that the density of negative integers is $1/2$.
The negative primes are exactly those with $q\equiv2\pmod3$, which also account for half the primes by the prime number theorem in arithmetic progressions \cite{Davenport}.

If $a+b=3^k$, write $a=3^t u,b=3^s v$ with $3\nmid uv$. We must have $t=s<k$, since otherwise dividing by the smaller power gives a contradiction modulo $3$.
Thus $u+v\equiv0\pmod3$, and hence
\[
f(a)f(b)=(c_t)^2\chi_3(u)\chi_3(v)=\chi_3(u)\chi_3(v)=-1.
\]
Therefore for all $k\ge1$ we have $R_{--}(3^k)=R_{++}(3^k)=0$.
\end{proof}

\subsubsection{A universal sufficient condition above density one half}
If the set $A$ of negative integers has natural density $d>1/2$, then without any multiplicativity assumption,
\[
R_{--}(N)\ge2|A\cap[1,N-1]|-(N-1)
=(2d-1)N+o(N)>0.
\]
Thus density strictly greater than one half guarantees existence, while equality does not give a universal guarantee.

\subsection{Degenerate probabilities and the prime-power rule}
At a deterministic endpoint, assigning the value $-1$ to every prime still leaves a choice of values at higher prime powers.
When every prime takes the value $-1$, there is no randomness at the primes, but the prime-power rule remains to be specified.
Completely multiplicative extension, extension by zero off the squarefree integers, and assigning $-1$ at every higher prime power give three different functions.
Each deterministic case must be studied with its own prime-power rule. The fair-model estimates apply to independent fair signs.
If zero values are allowed, the four counts of nonzero sign pairs sum to the number of pairs in the support, which determines the center.

For additive questions about the deterministic unitary function $g_0(n)=(-1)^{\omega(n)}$, see Tian and Yang \cite{TianYang}.
This function has the convolution representation $g_0=\mu*\kappa$, where $\kappa(1)=1$,
$\kappa(q^j)=1-j$, and for every $s>1/2$
$\sum_d|\kappa(d)|d^{-s}<\infty$. See Lemma 2.3 of their paper.
Under the Riemann hypothesis, the classical estimate $M_\mu(x)=O_\varepsilon(x^{1/2+\varepsilon})$ therefore gives
\[
\left|\sum_{n\le x}g_0(n)\right|
\le\sum_{d\le x}|\kappa(d)|\,|M_\mu(x/d)|
\ll_\varepsilon x^{1/2+\varepsilon}.
\]

\subsection{Open problems}\label{sec:questions}
The results leave the following questions.
\begin{question}[One-sided bias]
If $0\le\rho_q\le1/2$ and $\sum_qm_q/q=\infty$, does the quarter law hold almost surely without requiring $A_N\to\infty$? In particular, what happens if $\rho_q\to0$?
\end{question}

\begin{question}[Scale condition]
Is there a verifiable condition weaker than $A_N\to\infty$ that ensures $C_f(N)=o(N)$ along all integers?
\end{question}

\begin{question}[An inverse problem for reflection correlation]
If $|C_f(N)|\gg N$, must $f$ be close to a character at scale $N$, or have some other quantifiable multiplicative structure?
\end{question}

\appendix
\noindent Appendix~\ref{app:cm-technical} proves two technical lemmas used in Section~\ref{sec:cm-proof}.
Appendix~\ref{app:biased} proves the fixed-bias estimates for completely multiplicative functions in Corollary~\ref{cor:fixed}.
Appendix~\ref{sec:blocks} gives the block construction for Theorem~\ref{thm:counter}.

\section{Two technical lemmas for the completely multiplicative extremes}\label{app:cm-technical}
This appendix proves the local prefix-weight bound and the relative moderate-deviation estimate.

\subsection{Local weight of a prefix}
We use the notation $Q_N(a)$, $L,l,U,D$ and the groups from Section~\ref{sec:cm-upper}.
\begin{proof}
Write $u=t/g$. A squarefree $d\mid u$ is coprime to $g$. If $d\le N/g$,
the number of positions satisfying $gd\mid Q_N(a)$ is at most $2(N/gd)2^{\omega(g)+\omega(d)}$.
Below $2^{\omega(g)}=\exp(o(D))$, and $M_t^2\le2M_t$.
First consider $u>N/g$. Multiply the prime factors of $u$ in increasing order, taking the longest prefix $d$ that does not exceed $N/g$.
Let $w$ be the next prime. Set $k=\lfloor L/\log w\rfloor$ and $Y=e^{L/k}$.
Then $d>N/(gY)$, all prime factors of $d$ are $\le Y$, and $u/d<Nw$, so at most $k+2$ prime factors remain.
Group by $k$ and insert the factor $(dgY/N)^\sigma\ge1$ into the upper-bound sum over $d$, with
$\sigma=s/\log Y$. The preceding congruence count and the Euler product over squarefree divisors bound the contribution from one group by
\[
 \frac{N}{g}\exp(o(D))\, v^{k+2}e^{-(k-1)s}g^\sigma
       \exp(2vJ+2vH_s),\quad
 J=\sum_{p\le N}p^{-1},\quad H_s=\sum_{p\le Y}(p^\sigma-1)/p .
\]
Here $J=O(l)$. Take $\beta=1$ in the unrestricted case and $\beta=e^3$ in the smooth case.
When $\log Y\ge(\log(2\beta v))^2$, take $s=\log(\beta v)$.
Partial summation using Mertens' first theorem \eqref{eq:mertens-weighted} gives
\[
 H_s=\int_0^s (e^t-1)\,\frac{dt}{t}+O(e^s/\log Y)
       =(1+o(1))e^s/s .
\]
The error is uniform because $(e^{\sigma x}-1)/x$ is monotone and $\sigma=o(1)$.
Thus $2vH_s\le(6\alpha_0^2 \beta+o(1))D$, and the prefix factor is $O_\beta(v^3)\beta^{-k}$.
In the unrestricted case retain $g^{o(1)}$. In the smooth case this factor is $\exp(o(D))$,
and $k\ge D/0.30-1$, so $\beta^{-k}\ll e^{-10D}$.

When $\log Y<(\log(2\beta v))^2$, take $s=\tfrac32\log(2\beta v)$ instead.
Convexity gives $H_s\ll e^s$. Since $k\gg L/l^2$, we have $2vH_s=O(v^{5/2})=o(k)$.
Also $g^\sigma\le e^{ks/4}$, so $v^ke^{-ks}g^\sigma\le e^{-k\log v/8}$.
The decay in this range is stronger than $e^{-10D}$. There are at most $O(L)$ possible values of $k$.
In the unrestricted case with $u\le N/g$, enumerate $d=u$ directly. The contribution is at most
$(N/g)\exp(o(D)+2vJ)$. In the smooth case, set $y=e^{0.30U}$ and
$\lambda=\log(\beta v)/\log y$. Using $M_t^2\le4$ and counting smooth $u$ gives at most
\[
 4(N/g)^{1-\lambda}\prod_{p\le y}(1+v p^{\lambda-1}).
\]
By the same partial summation, the logarithm of the product is $O(D)$, whereas
$\lambda\log(N/g)\sim D\log(\beta v)/0.30$, which absorbs that term.
\end{proof}

\subsection{Relative moderate deviations for one or two targets}
Use the standardized martingale and the good event $\mathcal G$ from Section~\ref{sec:cm-lower}, with the parameters of Lemma~\ref{lem:cm-relative}.
\begin{proof}
\par\smallskip\noindent\textbf{Step 1. Stopping and exponential tilting.}
If the next step has $|a_p|>X^{-\nu}$, or if including it would give $\sum|a_p|^2>4$, set that step and all subsequent steps to zero.
Denote the stopped predictable vectors by $b_p$. There is no stopping on $\mathcal G$. Set $w=z(1,\ldots,1)$,
$S=\sum b_p\xi_p$,  and $K_w=\sum\log\cosh(w\cdot b_p)$.
Iterated conditional expectations show that the density $Z_w=\exp(w\cdot S-K_w)$ has mean $1$.
Since $Z_w^2/Z_{2w}\le\exp(2\sum(w\cdot b_p)^2)\le e^{8|w|^2}$,
we have $\EE Z_w^2\le X^{32}$. Under $d\PP_w=Z_w\,d\PP$,
Cauchy--Schwarz and \eqref{eq:cm-good} give $\PP_w(\mathcal G^c)\ll_h X^{-24}$.
At each step the conditional mean is $m_p=\tanh(w\cdot b_p)$ and the conditional variance is $1-m_p^2$.
Write $M=\sum b_p(\xi_p-m_p)$, $D_w=\sum b_p m_p$, and
$\Gamma_w=\sum(1-m_p^2)b_p b_p^{\mathsf T}$. On $\mathcal G$, use
$\sum |b_p|^4\le4X^{-2\nu}$ and Taylor expansion to obtain
\begin{align*}
 \|\Gamma_w-I_d\|_{\rm op}&\ll X^{-\nu}+|w|^2 X^{-2\nu},\\
 |D_w-w|&\ll |w|X^{-\nu}+|w|^3 X^{-2\nu},\\
 |K_w-|w|^2/2|&\ll |w|^2 X^{-\nu}+|w|^4 X^{-2\nu}.
\end{align*}

\par\smallskip\noindent\textbf{Step 2. Covariance and Gaussian replacement.}
Under the tilted law, project $M$ onto any unit direction. The martingale differences $\Delta_p$ satisfy $|\Delta_p|\le2X^{-\nu}$,
and the total conditional variance $H$ is at most $4$ on every path. The estimates above give $\EE_w|H-1|\ll_h X^{-\nu}$.
Scale each martingale difference by a predictable factor $\alpha_p\in[0,1]$ so that the conditional variance at that step is
$v'_p=\min(v_p,1-H'_{p-})$, where $v_p=\EE_w(\Delta_p^2\mid\text{previous signs})$,
and $H'_{p-}$ is the sum of the previous truncated variances. This gives $H'=\min(H,1)$ at the end.
Add a Gaussian term $G\sqrt{1-H'}$ to the truncated sum, using an independent standard normal variable. By the martingale isometry and
$(1-\alpha_p)^2v_p\le v_p-v'_p$, the mean square difference between the resulting variable and the original projection is at most $\EE_w|H-1|$.
For $f(x)=e^{iux}$ write $P_tf(x)=\EE f(x+\sqrt t\,G)$.
In $P_{1-H'_p}f$ at step $p$, use a centered Gaussian of conditional variance $v'_p$ to replace
$\alpha_p\Delta_p$. Since $H'_p$ is predictable, the conditional expectations of the quadratic Taylor polynomials agree.
Both third absolute moments are $\ll X^{-\nu}v'_p$, so the error is $\ll |u|^3X^{-\nu}v'_p$.
Gaussian convolution recovers $P_{1-H'_{p-}}f$. Sum over the steps and include the preceding mean square error.
This bounds the difference between the characteristic function of the original projection and $e^{-u^2/2}$ by
$\ll_h |u|X^{-\nu/2}+|u|^3X^{-\nu}$, uniformly in the direction.

\par\smallskip\noindent\textbf{Step 3. Smoothing rectangles and comparing characteristic functions.}
Put $B=X^{\nu/40}$ and $b=1/B$. Truncate axis-parallel rectangles to $[-B,B]^d$.
By $\EE_w|M|^2\le4$ and the normal second moment, this costs $O(B^{-2})$. Smooth the indicator of the rectangle with a Gaussian kernel of coordinate standard deviation $b$.
Before smoothing, its Fourier transform has absolute value $\ll B^2$. At frequency $\theta$ take
$u=|\theta|$. Fourier inversion bounds the difference of the smoothed probabilities by
$\ll_h B^2(b^{-3}X^{-\nu/2}+b^{-5}X^{-\nu})$.
Bound the original rectangle by expanding or contracting every face by $\sqrt b$ and allowing for the probability $\ll e^{-1/(2b)}$ that the Gaussian perturbation exceeds this distance.
The standard normal probability of the boundary layer is $\ll\sqrt b$, which also covers rectangles too narrow to contract.
Thus the probabilities for $M$ and a standard normal vector differ by $O_h(X^{-\nu/80})$ on every rectangle.
On $\mathcal G$ we have $|D_w-w|\ll_h X^{-\nu/2}$. Moving the faces again and allowing for $\PP_w(\mathcal G^c)$
gives the same bound for $Y=S-w=M+D_w-w$.

\par\smallskip\noindent\textbf{Step 4. Undoing the tilt and bounding the relative error.}\par
Set $f_z(v)=e^{-z\sum_i v_i}\ind_{\{\min_i v_i\ge0\}}$. This is the average of the indicator of the rectangle
$\prod_i[0,t_i]$ over the independent exponential densities $\prod_i ze^{-zt_i}\,dt_i$.
Thus $\EE_w f_z(Y)=(e^{z^2/2}u_X)^d+O_h(X^{-\nu/80})$,
with the main term computed by integration against independent standard normal variables and of size $\asymp(1+z)^{-d}$.
Undoing the tilt on $\mathcal G$ contributes the factor
$e^{-|w|^2/2}e^{-w\cdot Y}(1+O_h(X^{-\nu/2}))$.
The stopping rule leaves the original sum unchanged on this event. Restoring the remaining paths adds an error $O_h(X^{-80})$ under the original law,
and costs $O_h(X^{-24})$ in the tilted integral because $0\le f_z\le1$.
Multiplying by $e^{-dz^2/2}$ and using $(1+z)^d X^{-\nu/80}\ll_h X^{-\nu/100}$ gives the relative error.
\end{proof}

\section{Errors in the completely multiplicative model with fixed bias}\label{app:biased}
We prove the fixed-bias completely multiplicative conclusions of Corollary~\ref{cor:fixed}. Fix $0<\rho<1$ and put
$\theta=1-2\rho$ and $m=\min(\rho,1-\rho)$. The upper bound for the reflection sum follows from \eqref{eq:weight-reflection}
and Mertens' second theorem \eqref{eq:mertens-reciprocal}:
\[
 |C_f(N)|\ll_{\rho,\mathcal X} N\left(\frac{\log\log N}{\log N}\right)^{4m}.
\]

\subsection{Marginal bias and the counting error}
\leavevmode\par
\medskip\noindent\textbf{Step 1. Prime averages and the one-point asymptotic.}
Now suppose $\rho\ne1/2$. By Hoeffding's inequality and
$\sum_{q\le x}(\log q)^2\le x(\log x)^2$,
\[
\PP\left(\left|\sum_{q\le x}(\xi_q-\theta)\log q\right|>x^{3/4}\right)
\le2\exp\left(-\frac{c\sqrt x}{(\log x)^2}\right).
\]
Summing over integer $x$ and combining with the classical error term in the prime number theorem gives, almost surely for each fixed $A>0$,
\[
\sum_{q\le x}f(q)\log q=\theta x+O_{A,f}(x/(\log x)^A).
\]
Take a countable intersection for integer $A$. Together with $|f|\le1$, the mean-value theorem for multiplicative functions of Granville and Koukoulopoulos
\cite[Theorem~1]{GranvilleKoukoulopoulos} applies and gives
\[
S_f(x)
=c_{\rho,f}\frac{x}{(\log x)^{2\rho}}
+O_{\rho,f}\left(\frac{x}{(\log x)^{2\rho+1}}\right),
\qquad c_{\rho,f}=\frac{H_f}{\Gamma(1-2\rho)},
\]
where
\[
H_f=\prod_q(1-1/q)^\theta(1-f(q)/q)^{-1}>0.
\]
The product is taken in increasing prime order. Its logarithm is the convergent centered independent series
$\sum_q(f(q)-\theta)/q$ plus an absolutely convergent remainder, so almost surely $0<H_f<\infty$.
Thus $c_{\rho,f}$ is positive when $\rho<1/2$ and negative when $\rho>1/2$.
Aymone \cite{AymoneSD} studied the logarithmic scale of one-point sums in the squarefree model with fixed bias. The completely multiplicative result here uses the mean-value theorem just applied.

\medskip\noindent\textbf{Step 2. Comparing the two scales via the counting identity.}
Substitute the one-point asymptotic and the reflection bound into \eqref{eq:identity}. When $0<\rho<2/3$ and $\rho\ne1/2$,
we have $4m>2\rho$, so $C_f(N)=o(N/(\log N)^{2\rho})$. The leading errors in the two same-sign counts are
$\pm S_f(N-1)/2$. The opposite-sign errors have smaller order, giving the asymptotic equivalence in the corollary.
When $2/3\le\rho<1$, the lower bound follows from $R^f_{++}-R^f_{--}=S_f(N-1)$.
The upper bound follows from $4m\le2\rho$ and the same counting identity.
\begin{remark}[The fair point]
When $\rho=1/2$, the required result follows from Section~\ref{sec:fair-proof}.
The reciprocal gamma function vanishes at $0$. The asymptotic equivalence with a nonzero leading term therefore applies only for fixed $\rho\ne1/2$.
\end{remark}

\section{A block counterexample with divergent total probability weight}\label{sec:blocks}
We prove Theorem~\ref{thm:counter}. Alternate blocks of fair prime signs with character-matching blocks
so that in every stage, each realization has a target with no same-sign representations. Modifying the prime-power rule at $2$ gives
the non-completely multiplicative version. We then modify the probabilities to lie strictly between zero and one and verify the zero mean and the divergence for fixed real characters.

\begin{lemma}[Character matching for a finite prime configuration]\label{lem:match}
Given arbitrary signs $s_q$ at the finitely many odd primes $q\le H$, there exist arbitrarily large primes $P\equiv7\pmod8$ such that
\[
(2/P)=1,\qquad(q/P)=s_q\quad(q\le H,\ q\text{ is an odd prime}).
\]
\end{lemma}
\begin{proof}
Since $P\equiv3\pmod4$, quadratic reciprocity gives
$(q/P)=(P/q)(-1)^{(q-1)/2}$.
For each $q$ choose a nonzero class $a_q$ with $(a_q/q)=s_q(-1)^{(q-1)/2}$. Use the Chinese remainder theorem to combine these classes with $P\equiv7\pmod8$ into a reduced arithmetic progression.
Dirichlet's theorem on primes in arithmetic progressions gives arbitrarily large primes in it. The second supplementary law of quadratic reciprocity gives $(2/P)=1$
\cite[Sections~7.9 and 9.8]{Apostol}.
The modulus is fixed at each application of Dirichlet's theorem, which gives the required existence.
\end{proof}
\begin{remark}[Range of the construction parameters]
The recursion below requires only arbitrarily large matching primes and allows the stage parameters to grow rapidly. It gives no effective upper bound on the size of a target prime.
\end{remark}

\subsection{Fixing the distribution in advance: the stage recursion}
Set $\xi_2=1$ and take $L_0=2$ as the initial upper bound. The procedure enumerates finite configurations and fixes all probabilities before sampling.
The input at stage $j$ is $L_{j-1}$ together with the probability rules already assigned to the primes up to it.
\begin{enumerate}
\item \textbf{Choose a fair block.} Choose a finite $T_j>L_{j-1}$ such that
$\sum_{L_{j-1}<q\le T_j}q^{-1}\ge1$. On
$B_j=\{q:L_{j-1}<q\le T_j\}$ assign independent fair signs and set the current upper bound to $U=T_j$.
\item \textbf{Enumerate finite configurations.} List all configurations on the odd primes $q\le T_j$,
$\boldsymbol{\epsilon}=(\epsilon_q)\in\{\pm1\}^{\{q\le T_j:q\text{ is an odd prime}\}}$.
For each configuration, specify the signs $s_q$ of all odd primes $q\le U$ up to the current bound $U$.
For $q\le T_j$ take $s_q=\epsilon_q$. For $T_j<q\le U$ use the signs fixed earlier in this stage.
The sign at $2$ is always $1$ and is matched by Lemma~\ref{lem:match} as well.
\item \textbf{Extend the matching block.} Apply Lemma~\ref{lem:match} to choose $P>U+1$, requiring also
$(\log P)^3>T_j$. At every as yet unassigned prime $U<q<P$, set $\xi_q=\chi_P(q)$ deterministically,
where $\chi_P$ is the quadratic character modulo $P$. Update $U=P-1$ and continue to the next configuration.
\item \textbf{Proceed to the next stage.} After processing all configurations, set $L_j=U$ and proceed.
\end{enumerate}
The primes are thus divided in succession into fair blocks $(L_{j-1},T_j]$ and character-matching intervals.
Once set, the rule for each block remains fixed.
Each stage takes finitely many steps and the upper bounds grow without limit, so every prime is eventually assigned a rule.
The target prime $P$ itself is assigned its rule later, which has no effect on values at $n<P$.
The procedure uses only the enumerated finite configurations. It never reads the outcome of sampling, so the distribution is fixed in advance.
The prime signs are independent, with the deterministic coordinates viewed as degenerate independent variables. Since the fair blocks are disjoint,
\begin{equation}\label{eq:total}
\sum_qm_q/q=\frac12\sum_j\sum_{q\in B_j}1/q=\infty.
\end{equation}

\subsection{Infinitely many failing targets in every realization}
Fix any realization. At each stage its actual configuration on the primes $q\le T_j$ is enumerated.
At the corresponding target $P$, all realized signs at $q<P$ equal $\chi_P(q)$. The listed configuration matches the earlier blocks, the matching lemma covers the deterministic intervals already assigned in this stage, and the new interval is assigned by definition. The sign at $2$ also matches.
The completely multiplicative extension $H$ thus satisfies $H(n)=\chi_P(n)$ for all $1\le n<P$.
Hence $H(a)H(P-a)=\chi_P(-1)=-1$, so there are no same-sign decompositions.
By reflection symmetry, the two opposite-sign counts each account for half the pairs and both same-sign counts are zero.
Each stage provides at least one such target, and the stage bounds increase strictly. Thus there are infinitely many failing targets.
\begin{remark}[The cutoff scale at a failing target]
For these $P$, all primes assigned fair signs are either at most $T_j$ or at least $P$, while $(\log P)^3>T_j$.
Thus at most the prime $P$ itself contributes to the scale sum, and $0\le A_P=m_P/P\le1/(2P)$.
If the target is the last prime selected in its stage, $P=L_j+1$ can belong to the next fair block, so this endpoint contribution must be retained.
These targets are therefore outside the scale hypothesis of Theorem~\ref{thm:main}.
\end{remark}

\subsection{A genuinely non-completely multiplicative version}
Define
\[
f(n)=H(n)h(n),\qquad h(n)=\begin{cases}-1,&v_2(n)=2,\\1,&v_2(n)\ne2.\end{cases}
\]
Then $f(2)=1,f(4)=-1$. Multiplicativity on coprime arguments and the prime sign probabilities are preserved.
At a matched target $P=8m+7$, a pair changes its product sign only if exactly one summand has $v_2=2$.
These positions are $a\equiv4\pmod8$ or $a\equiv3\pmod8$. There are $2(m+1)=(P+1)/4$ of them in $1\le a<P$.
Therefore
\begin{equation}\label{eq:largeC}
C_f(P)=-(P-1)+2(P+1)/4=-(P-3)/2.
\end{equation}
By \eqref{eq:identity}, the four quarter-law asymptotics cannot all hold.

\subsection{Making all probabilities nondegenerate}
Independently flip each originally deterministic sign with probability $\delta_q>0$, chosen so that
\[
\sum_q\delta_q<\kappa.
\]
Leave the fair blocks unchanged. All probabilities are now strictly between $0$ and $1$, and the total weight still diverges.
By the union bound, the probability that none of the original deterministic signs flips is at least $1-\kappa$.
On this event there are still infinitely many failures.
\begin{remark}[Quantifiers for the modified probabilities]
After the modification, infinitely many failures occur on an event of probability at least $1-\kappa$. This suffices to disprove the general assertion of almost-sure success.
The original construction gives failure for every realization, whereas the modification gives this lower bound on its probability.
\end{remark}

\subsection{Zero mean and divergence for fixed real characters}
Return to the original construction without the additional flips, and let $B=\bigcup_j B_j$.
For every fixed real Dirichlet character $\chi$, the centered series of independent terms
\[
\sum_{q\in B}\frac{\xi_q\chi(q)}q
\]
has summable variances and therefore converges almost surely, by the standard convergence theorem for independent series or the $L^2$-bounded martingale convergence theorem \cite{Durrett}.
On the other hand $\sum_{q\in B}1/q=\infty$, so
\begin{equation}\label{eq:fixedchar}
\sum_q\frac{1-f(q)\chi(q)}q=\infty
\end{equation}
holds almost surely, since all terms are nonnegative and it suffices to look at $B$. The real characters are countable, so a single probability-$1$ event covers all fixed real characters.

For the principal character modulo $1$, equation \eqref{eq:fixedchar} gives $\sum_q(1-f(q))/q=\infty$,
and $H(q)=f(q)$ gives the same divergence condition for $H$. In its real-valued form, Wirsing's theorem
\cite{Wirsing} (see also \cite[Theorem A, p.~4]{Pollack}) states that a multiplicative function
$u:\mathbb N\to[-1,1]$ satisfying $\sum_q(1-u(q))/q=\infty$ has $S_u(N)=o(N)$.
Here $H$ is completely multiplicative, and $f=Hh$ is multiplicative on coprime arguments and has values $\pm1$ at every prime power.
Both functions therefore satisfy the theorem.
Thus almost surely the positive and negative integers each have density one half, while the reflection correlations remain large along the moving targets.
For the non-completely multiplicative model, \eqref{eq:identity} gives along the sequence of matched primes
\[
R_{--}^f(P)=P/8+o(P),\quad R_{++}^f(P)=P/8+o(P),\qquad
R_{-+}^f(P)=R_{+-}^f(P)=3P/8+O(1).
\]
\begin{remark}[Fixed characters and moving structure]
By \eqref{eq:fixedchar} the weighted distance from each fixed real Dirichlet character diverges.
Characters varying with the target, complex characters, and twists by $n^{it}$ form a wider class. Uniform control over this class and an inverse description require further work.
\end{remark}

\section*{Use of AI tools}
The mathematical exploration and proof drafts relied primarily on the AI tool ChatGPT6-Astra. ChatGPT6-sol and DeepSeek were also used in preparing the manuscript, for drafting, translation, and preliminary literature searches.

\bibliographystyle{amsplain}
\providecommand{\bysame}{\leavevmode\hbox to3em{\hrulefill}\thinspace}
\providecommand{\MR}{\relax\ifhmode\unskip\space\fi MR }
\providecommand{\MRhref}[2]{%
  \href{http://www.ams.org/mathscinet-getitem?mr=#1}{#2}
}
\providecommand{\href}[2]{#2}

\end{document}